\documentclass[12pt]{amsart}

\usepackage{amscd, amssymb, tikz,enumerate, xcolor,textcomp,ragged2e,blindtext,multicol,mathrsfs}
\usetikzlibrary{calc,3d,decorations.markings}
\usetikzlibrary{arrows}
\usepackage[colorlinks=true, citecolor=blue!50!black, linkcolor=red!50!black]{hyperref}

\usepackage{geometry}
\pdfpagewidth=\paperwidth
\pdfpageheight=\paperheight

\newtheorem{thm}{Theorem}[section]
\newtheorem{cor}[thm]{Corollary}
\newtheorem{lem}[thm]{Lemma}
\newtheorem{prop}[thm]{Proposition}

\theoremstyle{definition}
\newtheorem{dfn}[thm]{Definition}

\theoremstyle{remark}
\newtheorem{rmk}[thm]{Remark}
\newtheorem{rmks}[thm]{Remarks}
\newtheorem{example}[thm]{Example}
\newtheorem{examples}[thm]{Examples}

\newtheorem{nota}[thm]{Notation}

\numberwithin{equation}{section}

\newcommand{\ZZ}{\mathbb{Z}}
\newcommand{\NN}{\mathbb{N}}

\newcommand{\FF}{\mathbb{F}}

\newcommand{\Gg}{\mathcal{G}}

\newcommand{\Tt}{\mathcal{T}}

\newcommand{\id}{\operatorname{id}}

\newcommand{\Aut}{\operatorname{Aut}}

\begin{document}

\title[Embeddability of higher-rank graphs]{Embeddability of higher-rank graphs in groupoids, and the structure of their $C^*$-algebras}

\author{Nathan Brownlowe}
\address{Nathan Brownlowe\\ School of Mathematics and Statistics \\
The University of Sydney\\
NSW 2006 \\
Australia.}
\email[Brownlowe]{nathan.brownlowe@sydney.edu.au}

\author{Alex Kumjian}
\address{Alex Kumjian\\ Department of Mathematics (084)\\ University
of Nevada\\ Reno NV 89557-0084\\ USA}
\email[Kumjian]{alex@unr.edu}

\author{David Pask}
\address{David Pask\\ School of Mathematics \&
Applied Statistics  \\
University of Wollongong\\
NSW  2522\\
Australia}
\email[Pask]{david.a.pask@gmail.com}
\author{Aidan Sims}
\address{Aidan Sims\\ School of Mathematics \&
Applied Statistics  \\
University of Wollongong\\
NSW  2522\\
Australia}
\email[Sims]{asims@uow.edu.au}

\subjclass[2020]{46L05 (primary); 18D99 (secondary)}
\keywords{$k$-graph; $C^*$-algebra; fundamental groupoid; fundamental group; $\tilde{A_2}$-buildings}

\thanks{
This research was supported by the Australian Research Council Discovery Project DP220101631.
A.K.\ was partially supported  by  Simons Foundation Collaboration grant 353626.}
\date{\today}

\dedicatory{This paper is dedicated to our friend and mentor Iain Raeburn, whose guidance had a profound impact on us all. We'd love to say more, but he'd have been squirming enough already.}

\begin{abstract}
We show that the $C^*$-algebra of a row-finite source-free $k$-graph is
Rieffel--Morita equivalent to a crossed product of an AF algebra by the
fundamental group of the $k$-graph. When the $k$-graph embeds
in its fundamental groupoid, this AF algebra is a Fell algebra; and
simple-connectedness of a certain sub-1-graph characterises when this Fell
algebra is Rieffel--Morita equivalent to a commutative $C^*$-algebra. We
provide a substantial suite of results for determining if a given $k$-graph
embeds in its fundamental groupoid, and provide a large class of examples,
arising via work of Cartwright, Robertson, Steger  et al.\ from the theory
of $\tilde{A_2}$-groups, that do embed.
\end{abstract}

\maketitle

\section{Introduction}

Since their introduction \cite{KP2} higher-rank graphs, or $k$-graphs, have been a source of interesting new higher-dimensional phenomena: in algebra \cite{APCaHR,CFaH,Ros}, dynamics \cite{PRW,Kak,SZ,Sp}, $C^*$-algebras \cite{AB,CaHS,RSS}, $K$-theory \cite{E,G,PRenS}, topology \cite{PQR1,PQR2,KKQS,KPSW} and geometry \cite{RobSteg1,V,KV}. However, many natural questions about their structure theory remain difficult to unravel.

One such question, and the primary motivation for this paper, is: when can a $k$-graph $C^*$-algebra be realised, modulo Rieffel--Morita equivalences, as a crossed product of a commutative $C^*$-algebra? For $1$-graphs, the answer is ``always:" given a row-finite source-free directed graph $E$, the middle two authors showed \cite{KP1} that the $C^*$-algebra of its universal cover $F$ is Rieffel--Morita equivalent to a commutative AF algebra, and there is an action of the fundamental group $\pi_1(E, v)$ on $C^*(F)$ whose crossed product is Rieffel--Morita equivalent to $C^*(E)$. For $k$-graphs, the answer is more nuanced, and is related to two other intriguing structural questions: when does a $k$-graph embed in its fundamental groupoid, and when is the boundary of its universal cover Hausdorff?

Our main $C^*$-algebraic theorem, Theorem~\ref{thm:mainC*result}, clarifies the relationships between these questions: the $C^*$-algebra $C^*(\Lambda)$ of any row-finite source-free $k$-graph is a crossed product of an $AF$ algebra $C^*(\Sigma)$ by the fundamental group of $\Lambda$; if $\Lambda$ embeds in its fundamental groupoid, then the AF algebra $C^*(\Sigma)$ is a Fell algebra; and if, additionally, a naturally-arising sub-1-graph of $\Sigma$ is simply connected, then the boundary of $\Sigma$ is Hausdorff, and $C^*(\Sigma)$ is Rieffel--Morita equivalent to a commutative AF algebra. The point is that the first part of the program of \cite{KP1} above goes through smoothly for $k$-graphs: every $k$-graph $\Lambda$ has a fundamental group $\pi(\Lambda)$ \cite{PQR1} and a universal cover $\Sigma$ \cite{PQR2} that carries an action of $\pi(\Lambda)$, and when $\Lambda$ is row-finite and source-free, the resulting crossed product is Rieffel--Morita equivalent to $C^*(\Lambda)$ \cite{KP2}. Our main contribution is the analysis of $C^*(\Sigma)$.

Motivated by this, we study the question of when a $k$-graph $\Lambda$ embeds in its fundamental groupoid.  Many $k$-graphs do not embed: we give three examples in Section~\ref{subsec:non-embeddings}; and any $k$-graph containing a copy of one of these (of which there are many) also fail to embed. So we focus on checkable sufficient conditions. We show that \emph{singly connected} $k$-graphs always embed (Proposition~\ref{prop:1-cnc-emb}(ii)), and highlight a surprising difference between $k$-graphs and $1$-graphs: universal covers of $k$-graphs need not be singly connected. We include a proof that $1$-graphs always embed (Theorem~\ref{thm:1gcase}). We then show that many standard $k$-graph constructions preserve embeddability: coverings (Proposition~\ref{prop:lift}), affine pullbacks, cartesian products, crossed-products and skew-products (Corollary~\ref{cor:specialcases}), and action graphs (Corollary~\ref{thm:genaction}). The workhorse in this is Theorem~\ref{prop:prec} which exploits the universal property of the fundamental groupoid and fundamental group. In Proposition~\ref{prop:any vertex} we reduce the embeddability of a connected $k$-graph to group-embeddability of the subsemigroup based at any vertex. Using Dilian Yang's work \cite{Y4} on $k$-graphs and Yang--Baxter solutions, we show that there are many embeddable $k$-graphs for all $k$ (Lemma~\ref{lem:ybex}).  We are far from a complete answer to the embeddability question. Johnstone's general results \cite{J} characterise groupoid-embeddability of categories, but the hypotheses seem uncheckable: we gleaned no practical conditions---either necessary or sufficient---from Johnstone's work, beyond the neat result of Lawson--Vdovina \cite[Theorem~11.14]{LV} presented in Remark~\ref{rmk:LV}. One might hope for help from Ore's theorem \cite[Proposition II.3.11]{Dehornoy}, but by the factorisation property, no interesting $k$-graphs are Ore. Remarks \ref{rmk:confused}~and~\ref{rmk:more confused} indicate how much we still do not know about embeddability.

Finally, as definitive general results about embeddability are still beyond reach, we present a class of examples arising from the combinatorial objects used by Roberston and Steger to construct higher-rank Cuntz--Krieger $C^*$-algebras in \cite{RobSteg1} which first inspired the middle two authors to develop the concept of a higher-rank graph. We show in Theorem~\ref{thm:lambdaT} and Proposition~\ref{prop:skew-cover} that every $\tilde{A_2}$-group $\Gamma_\Tt$ yields $2$-graphs $\Lambda_\Tt$ and $\Sigma_\Tt$, the latter being a cover of the former, and in Corollary~\ref{cor:LambdaT embeds} that $\Lambda_\Tt$ embeds in its fundamental groupoid---Proposition~\ref{prop:lift} then shows that $\Sigma_\Tt$ embeds as well. We also prove that $\Sigma_\Tt$ is singly connected, and deduce that its $C^*$-algebra is type~I$_0$. The construction of $\Lambda_\Tt$ is related to a number of existing constructions. It is directly inspired by \cite[pp.135--136]{RobSteg1}. As discussed in \cite{CMSZ1, RobSteg1} a thick $\tilde{A_2}$-building $\mathscr{B}$ carrying a vertex-transitive action of a $\tilde{A_2}$-group $\Gamma_\Tt$ arises from a finite projective plane $(P, L)$, a bijection between $P$  and $L$, and a compatible triangle presentation $\Tt$ on $P$, the points of the projective plane, arising from the local structure of the building (see \cite[\S{3}]{CMSZ1}). The $\tilde{A_2}$-group $\Gamma_\Tt$ is generated by a set indexed by $P$ subject to the relations encoded in  $\Tt$. The $\tilde{A_2}$-building $\mathscr{B}$ is constructed as an augmented Cayley graph of $\Gamma_\Tt$ with 2-simplices given by $\Tt$. Our $\Lambda_\Tt$ is isomorphic to the $2$-graph obtained from \cite[Examples~1.7(iv)]{KP2} from the $0$--$1$ matrices $M_i$ of \cite[p.135]{RobSteg1} (see Remark~\ref{rmk:us vs RobSteg}). Geometric considerations suggest both that $\Sigma_\Tt$ should be simply connected, and therefore equal to the universal cover of $\Lambda_\Tt$, and that its topological realisation should coincide with that of $\mathscr{B}$, so it should have Hausdorff boundary; we leave this for future work. Our construction is also related to the construction of $k$-graphs from groups in \cite{MRV}, but cannot be recovered from it: the covering $2$-graphs in \cite{MRV} are products of trees rather than $\tilde{A_2}$-buildings.

\section{Background and preliminary results}\label{sec:background}

\subsection{Higher rank graphs}

We write $\NN$ for the additive monoid $\{ 0, 1 ,\ldots \}$. We denote the standard generators of $\NN^k \subset \ZZ^k$ by $\varepsilon_1, \dots ,\varepsilon_k$, and we write $n_i$ for the
$i^{\textrm{th}}$ coordinate of $n \in \NN^k$. We write $\mathbf{1}_k$ or just $\mathbf{1}$ for $(1, \ldots , 1) \in \NN^k$.

A \emph{$k$-graph} is a small category $\Lambda$ equipped with a functor $d :\Lambda\to\NN^k$ satisfying the \emph{factorization property}: whenever $d(\lambda) = m+n$, there exist unique $\mu,\nu\in\Lambda$ such that $d(\mu)=m$, $d(\nu)=n$, and $\lambda=\mu\nu$. This implies that $\Lambda$ is cancellative. We write $\Lambda^n := d^{-1} (n)$ for $n \in \NN^k$. When $d(\lambda)=n$ we say $\lambda$ has \emph{degree} $n$. The factorisation property implies that $\Lambda^0$ is the set of identity morphisms, which we call \emph{vertices}. Elements of $\bigcup_i \Lambda^{\varepsilon_i}$ are called \emph{edges}. For $u,v\in\Lambda^0$ we write $u\Lambda := r^{-1}(u)$, $\Lambda v:=s^{-1}(v)$, and $u\Lambda v:=u\Lambda \cap \Lambda v$.

\begin{nota}
For $\lambda \in \Lambda$ and $0 \le m \le n \le d(\lambda)$, we write $\lambda(m,n)$ for the unique element of $\Lambda$ such that $\lambda \in \Lambda^m \lambda(m,n)\Lambda^{d(\lambda) -n}$.
We define $\lambda(n) := \lambda(n,n) = s(\lambda(0,n))$.
\end{nota}

\begin{dfn} \label{dfn:connected}
The $k$-graph  $\Lambda$ is \emph{connected} if the equivalence relation $\sim$ on $\Lambda^0$ generated by $\{(u,v)\mid u\Lambda v\ne\emptyset\}$ is $\Lambda^0\times\Lambda^0$. A $k$-graph is \textit{strongly connected} if $u \Lambda v \neq \emptyset$ for all $u,v \in \Lambda^0$.
\end{dfn}

A \emph{morphism} $\phi : \Omega \to \Lambda$  between $k$-graphs is a functor such that $d_\Lambda ( \phi ( \lambda ) ) = d_\Omega ( \lambda )$ for all $\lambda \in \Omega$. A \emph{quasi-morphism} from a $k$-graph $( \Omega , d_\Omega )$ to an $\ell$-graph $( \Lambda , d_\Lambda )$ is a pair $( \phi , f)$ consisting of a functor $\phi : \Omega \to \Lambda$ and a homomorphism $f: \NN^k \to \NN^\ell$ such that $d_\Lambda \circ \phi = f \circ d_\Omega$. If $\Lambda$ is a $k$-graph, then $\Lambda^{\NN{\bf 1}} := \{\lambda \in \Lambda : d(\lambda) \in \NN{\bf 1}\}$ is a $1$-graph and the natural inclusion $\Lambda^{\NN{\bf 1}} \hookrightarrow \Lambda$ together with the map $f: \NN \to \NN^k$ given by $f(n) := n{\bf 1}$ is a quasimorphism.

\begin{examples} \label{ex:bnex}
\begin{enumerate}[(i)]
\item Let $B_n$ be the directed graph with $B_n^{0} = \{u\}$, and $B^1_n = \{f_1 , \ldots , f_n\}$. Its path category $B_n^*$ is a $1$-graph, and coincides with the free semigroup $\mathbb{F}_n^+$ on $n$ generators.
\item Let $\Delta_k = \{ (m,n) \in \mathbb{Z}^k \times \mathbb{Z}^k : m \le n \}$. Define $r,s : \Delta_k \to
\operatorname{Obj} \Delta_k$ by $r (m,n) = m$,  $s (m,n) = n$,
and for $m \le m \le p \in \ZZ^k$ define $(m,n)(n,p)=(m,p)$ and $d (m,n) = n-m$. Then $( \Delta_k , d )$ is a $k$-graph where $\operatorname{Obj} \Delta_k$ is identified with $\{ (m,m) : m \in \mathbb{Z}^k \}
\subset \operatorname{Mor} \Delta_k$.
\item Similarly, $\Omega_k = \{(m,n) \in \mathbb{N}^k \times \mathbb{N}^k : m \le n\}$ is a sub-$k$-graph of $\Delta_k$.
\end{enumerate}
\end{examples}
%

\begin{example}[Skew-product graphs]\label{eg:skewproducts}
Let $\Lambda$ be a $k$-graph, $G$ a group, and $c: \Lambda \to G$ a $1$-cocycle (functor). Then the set $G \times_c \Lambda := \{ (g, \lambda) : g \in G, \lambda \in \Lambda \}$, under the structure maps
\begin{equation} \label{eq:skewprod}
s (g,  \lambda ) = ( g c ( \lambda ), s ( \lambda ) ), \quad r (g, \lambda ) = (g, r(\lambda )),  \quad
	(g, \lambda   ) \cdot (g c (\lambda), \mu) = (g, \lambda \mu ),  \quad d (g ,  \lambda ) = d ( \lambda ) .
\end{equation}
is a $k$-graph called the \emph{skew-product graph} \cite[Defintion 5.1]{KP2}. Left translation by $G$ on the first coordinate of $G \times_c \Lambda$ is an action of $G$ by $k$-graph automorphisms.

There are two equivalent conventions for skew-product graphs in the literature: the other is \cite[Defintion 6.3]{PQR2}. In \cite[Defintion 6.3]{PQR2} $\Lambda \times_c G := \{ (\lambda, g) : \lambda \in \Lambda, g \in G  \}$ with structure maps
\[
	s ( \lambda , g ) = ( s ( \lambda ) ,  g  ), \quad r ( \lambda , g ) = ( r(\lambda ) ,  c( \lambda)g  ),   \quad
	( \lambda ,  c( \lambda)g ) \cdot ( \mu , g ) = ( \lambda \mu , g ),  \quad d ( \lambda , g ) = d ( \lambda ) .
\]
It is simple to check that $\phi(g,  \lambda ) = (\lambda,  c( \lambda)^{-1}g^{-1})$ yields an isomorphism $\phi: G \times_c \Lambda \to \Lambda \times_c G$.
\end{example}

\begin{example}[Monoidal $2$-graphs]\label{eg:monoidalkgraphs}
The following class of 2-graphs was introduced in \cite[\S 6]{KP2} and later studied extensively by Yang et al.\ \cite{DPY, Y1, Y4}. Fix $n_1 , n_2 \ge 1$. Let $[n_i ]= \{ 1 , \ldots , n_i \}$, for $i=1,2$. Let $\theta : [n_1 ]\times [n_2] \to [n_2] \times [n_1]$ be a bijection. The monoidal $2$-graph $\mathbb{F}_\theta^+$ is the unique $2$-graph such that $(\mathbb{F}_\theta^+)^0 = \{v\}$, $(\mathbb{F}_\theta^+)^{\varepsilon_1} = \{ e_{1}, \ldots , e_{n_2} \}$, $(\mathbb{F}_\theta^+)^{\varepsilon_2} = \{ f_{1}, \ldots , f_{n_2} \}$, and
\begin{equation} \label{eq:thetadef}
		e_i f_j = f_{j'} e_{i'} \text{ whenever } \theta (i,j) = (j',i').
\end{equation}
\end{example}

\begin{rmk}
In their early papers Yang et al.\ define $\mathbb{F}_\theta^+$ in terms of a bijection $\theta : [n_1 ]\times [n_2] \to [n_1] \times [n_2]$ rather than $[n_1 ]\times [n_2] \to [n_2] \times [n_1]$.
\end{rmk}

An \emph{affine} map $f : \NN^\ell \to \NN^k$ is a map of the form $f(n) = An+p$ for $A \in M_{\ell,k} ( \NN )$ and $p \in \NN^\ell$. The following proposition unifies the pullback construction of \cite[Definition~1.9]{KP2} (case $p=0$) and the $p$-dual graph of \cite[Definition~3.2]{APS} (case $A=I$).

\begin{prop}[Affine pullbacks]\label{prop:gpb}
Let $(\Lambda,d)$ be a $k$-graph and let $f : \NN^\ell \to \NN^k$ be an affine map with $f(0)=p \in \NN^k$. Set $f^* (\Lambda ) = \{ ( \lambda , n ) : d ( \lambda )  = f(n) \} \subseteq \Lambda \times \NN^\ell$. Then $f^*(\Lambda)$ is an $\ell$-graph, with structure maps $r (\lambda , n ) = [ \lambda (0,p) , 0]$, $s  ( \lambda , n ) = [\lambda ( d(\lambda ) -p  , d(\lambda)),0]$,
\begin{equation} \label{eq:fstardef}
		(\lambda , m ) \circ  (\mu , n ) = ( \lambda ( 0 , d(\lambda)-p) \mu , m+n) \text{ if $s (\lambda , m ) = r ( \mu , n )$, }
\end{equation}
and $d_{f^*(\Lambda)} (\lambda,n) = n$. We have $f^*(\Lambda)^0 = \{ ( \lambda, 0 ) :  \lambda \in \Lambda^p \}$.
\end{prop}
\begin{proof}
As in \cite[Definition~1.9]{KP2}, the pullback $A^*(\Lambda)$ of $\Lambda$ by the homomorphism $A : \NN^\ell \to \NN^k$ is an $\ell$-graph. By \cite[Proposition~3.2]{APS}, its dual $p(A^*(\Lambda))$ is also an $\ell$-graph. As sets, $p(A^*(\Lambda)) = \{(\lambda, n) \in \Lambda \times \NN^\ell : d_{p\Lambda}(\lambda) = A n\} = \{(\lambda, n) \in \Lambda \times \NN^\ell : d(\lambda) = f(n)\} = f^*\Lambda$. Direct calculations show that this identification intertwines the structure maps above with those of $p(A^*(\Lambda))$.
\end{proof}

\begin{example}[Crossed-product graph] \label{ex:crossed product graph}
Let $\alpha : \ZZ^\ell \to \operatorname{Aut} \Lambda$ be an action of $\ZZ^\ell$ on a $k$-graph $\Lambda$. Then the set $\Lambda \times \NN^\ell$ with the structure maps
\[
r ( \lambda , m )  = ( r ( \lambda ) , 0 ) ,  \ s ( \lambda , m ) = ( \alpha^{-m} (s (\lambda ) ) , 0 ) , \
( \lambda , m )( \mu , n ) = ( \lambda \alpha^m ( \mu ) , m+n ),   \ d ( \lambda , m ) = ( d ( \lambda) , m )
\]
is a $(k + \ell)$-graph, called the \emph{crossed-product graph} $\Lambda \times_\alpha \NN^\ell$  (see \cite{FPS}).
\end{example}

\subsection{Fundamental groupoids, fundamental groups and universal covers}\label{sec:fun}

Every $k$-graph $\Lambda$ has a fundamental groupoid, defined as follows (see \cite[Section~19.1]{Schubert} or \cite[Section~3]{PQR1}).

\begin{dfn}\label{dfn:fundgpd}
Let $\Lambda$ be a $k$-graph. There exists a groupoid $\Pi(\Lambda)$ and a functor $i : \Lambda \to \Pi(\Lambda)$ such that $i(\Lambda^0)=\Pi(\Lambda)^0$, with the following universal property: for every functor $F$ from $\Lambda$ into a groupoid $\Gg$, there exists a unique groupoid homomorphism $\tilde{F} : \Pi(\Lambda) \to \Gg$ such that $\tilde{F} \circ i = F$. The pair $(\Pi(\Lambda), i)$ is unique up to canonical isomorphism, so we refer to any such groupoid $\Pi(\Lambda)$ as \emph{the fundamental groupoid} of $\Lambda$.
\end{dfn}

The assignment $\Lambda \mapsto \Pi (\Lambda)$ is a functor from $k$-graphs to groupoids. Note that $\Pi(\Lambda)$ is denoted $\Gg (\Lambda)$ in \cite{PQR1}, but this clashes with the notation for path groupoids in Section~\ref{sec:path gpd} and \cite{KP2}.

Each $k$-graph also has a fundamental group; the standard definition, for connected $k$-graphs, is as any one of the isotropy groups of its fundamental groupoid, as follows.

\begin{dfn}
Let $\Lambda$ be a $k$-graph. The \emph{pointed fundamental group} $\pi_1(\Lambda, v)$ of $\Lambda$ at $v \in\Lambda^0$ is the isotropy group $\pi_1 (\Lambda, v) := v \Pi ( \Lambda ) v$ of $\Pi(\Lambda)$ at $v$.
\end{dfn}

\begin{dfn}
For $X \not= \emptyset$, the \emph{pair groupoid} of $X$ is $T(X) := X \times X$, the simple transitive groupoid with unit space $\{(x,x) : x \in X\}$ identified with $X$; it has structure maps
\[
r(x, y) := x, \quad s(x, y) := y, \quad (x, y)(y, z) := (x, z), \quad (x, y)^{-1} := (y, x).
\]
\end{dfn}

\begin{rmk}\label{rmk:trans}
Suppose that $\Lambda$ is connected and let $v \in\Lambda^0$. Then there exists a function $w \mapsto \gamma_w$ from $\Lambda^0$ to $\Pi(\Lambda)_v$ such that $\gamma_v = v$ and $r(\gamma_w) = w$ for all $w$. Any such function $\gamma$ determines an isomorphism $(g, (u,w))\mapsto \gamma_u g \gamma_w$ from $\pi_1 (\Lambda, v)  \times T(\Lambda^0)$ to $\Pi(\Lambda)$ \cite[Corollary 6.5]{PQR2}.
\end{rmk}

The isomorphism class of $\pi_1(\Lambda, v)$ is independent of $v$ when $\Lambda$ is connected, but the isomorphisms $\pi_1(\Lambda, v) \to \pi_1(\Lambda, w)$ are noncanonical. However, each category also admits a canonical fundamental group (elsewhere called the universal group of the category) defined abstractly as the universal group generated by a cocycle on the category \cite{BKQ21, Higgins}; our next result shows how the two are related. In the next result, $\sim$ is the relation on $\Lambda^0$ from Definition~\ref{dfn:connected}.

\begin{prop}
Let $\Lambda$ be a $k$-graph. Let $V \subseteq \Lambda^0$ be a complete collection of representatives of the equivalence classes for $\sim$. Then the group
\[
\pi_1(\Lambda) := \bigoplus_{v \in V} \pi_1(\Lambda, v)
\]
has the following universal property: there is a $1$-cocycle $\kappa : \Lambda \to \pi_1(\Lambda)$ whose image generates $\pi_1(\Lambda)$ such that for any $1$-cocycle $c$ from $\Lambda$ to a group $G$ there is a unique homomorphism $\tilde{c} : \pi_1(\Lambda) \to G$ such that $c = \tilde{c} \circ \kappa$. If $G$ is a group and $c : \Lambda \to G$ a $1$-cocycle with the same universal property, then $\tilde{c}$ is an isomorphism of $\pi_1(\Lambda)$ onto $G$.
\end{prop}

\begin{proof}
First suppose that $\Lambda$ is connected, and fix $v \in \Lambda^0$. Let $\Pi(\Lambda)$ be its fundamental groupoid. Remark~\ref{rmk:trans} shows that there is an isomorphism $\varphi : \Pi (\Lambda )\to  \pi_1(\Lambda, v) \times_\gamma T(\Lambda^0)$ such that $\varphi(\gamma_u g \gamma_w^{-1}) = (g, (u,w))$ for all $(g, (u,w))$. Let $p_1 : \pi_1(\Lambda, v) \times_\gamma T(\Lambda^0) \to \pi_1(\Lambda, v)$ be the projection map $(\gamma, (u,w)) \mapsto \gamma$. Then there is a cocycle $\kappa : \Lambda \to \pi_1(\Lambda, v)$ given by $\kappa := p_1 \circ \varphi \circ i$. Fix $\eta := i(\lambda_1) i(\lambda_2)^{-1} \dots i(\lambda_n) \in \pi_1(\Lambda, v)$. Then since $\gamma_{r(\lambda_1)} = i(v) = \gamma_{s(\lambda_n)}$, we have
\[
\eta = \gamma_{r(\lambda_1)} i(\lambda_1)\gamma_{s(\lambda_1)} (\gamma_{r(\lambda_2)}^{-1} i(\lambda_2)\gamma_{s(\lambda_2)})^{-1} \cdots \gamma_{r(\lambda_n)}^{-1}i(\lambda_n)\gamma_{s(\lambda_n)}
= \varphi(\lambda_1)\varphi(\lambda_2)^{-1} \cdots \varphi(\lambda_n),
\]
so the range of $\kappa$ exhausts $\pi_1(\Lambda, v)$.
	
Let $c : \Lambda \to G$ be a $1$-cocycle into a group. The universal property of $\Pi(\Lambda)$ gives homomorphism $\rho : \Pi(\Lambda) \to G$ such that $\rho \circ i = c$. So $\rho \circ \varphi^{-1} : \pi_1(\Lambda, v) \times_\gamma T(\Lambda^0) \to G$ is a homomorphism. Define $\iota_v : \pi_1(\Lambda, v) \to \pi_1(\Lambda, v) \times_\gamma T(\Lambda^0)$ by $\iota_v(g) = (g, (v,v))$ and let $\tilde{c}:= \rho \circ \varphi^{-1} \circ \iota_v$. Since $\iota_v$ is inverse to $p_1$ on $\pi_1(\Lambda, v) \times \{(v,v)\}$, we have $\varphi^{-1} \circ \iota_v \circ p_1 \circ \varphi = \id_{\Pi(\Lambda)^v_v}$. So $\tilde{c}\circ \kappa(\lambda)) = \rho \circ i(\lambda) = c(\lambda)$.
	
Now suppose that $\Lambda$ is not connected. Let $[v]$ denote the equivalence class of $v \in \Lambda^0$ under $\sim$. For each connected component $\Lambda_v = [v]\Lambda[v]$, the above argument gives a $1$-cocycle $\kappa_v : \Lambda \to \pi_1(\Lambda, v)$. So $\bigoplus_{v \in V} \kappa_v : \Lambda \to \bigoplus_{v \in V} \pi_1(\Lambda, v)$ is a $1$-cocycle. Its image generates because its image in each summand generates. Any $1$-cocycle $c : \Lambda \to G$ into a group restricts to $1$-cocycles $c_v : \Lambda_v \to G$, which induce homomorphisms $\tilde{c}_v : \pi_1(\Lambda_v, v) \to G$ as above. Then $\tilde{c} := \bigoplus_v \tilde{c}_v$ satisfies $\tilde{c} \circ \kappa = c$.
	
For the final statement, observe that the universal property of $(G, \tilde{c})$ applied to $\kappa$ yields a homomorphism $\tilde\kappa : G \to \pi_1(\Lambda)$ that is inverse to $\tilde{c}$.
\end{proof}

The following definitions appear in \cite{PQR2}. We include them for completeness.

\begin{dfn}
	Let $\Lambda, \Sigma, \Gamma$ be $k$-graphs.
	\begin{enumerate}[(i)]
		\item A surjective $k$-graph morphism $p  : \Sigma \to \Lambda$ is a \emph{covering} if for all $v \in \Sigma^0$, $p$ restricts to bijections $\Sigma v \to\Lambda p(v)$ and $v \Sigma \to p(v) \Lambda$.
		\item A covering $p:\Sigma \to \Lambda$ is said to be \emph{connected} if $\Sigma$ (and hence $\Lambda$) is connected.
		\item If $p: \Sigma \to \Lambda$ and $q : \Gamma \to \Lambda$ are coverings, a \emph{morphism} from $(\Sigma,p)$ to $(\Gamma, q)$ is a $k$-graph morphism $\phi : \Sigma \to \Gamma$ such that $q \circ \phi = p$.
		\item A covering $p : \Sigma \to\Lambda$ is \emph{universal} if it is connected in the sense of~(ii), and for every connected covering
		$q  : \Gamma \to\Lambda$, there is a unique morphism $\phi : (\Sigma,p)\to (\Gamma,q)$ in the sense of~(iii).
	\end{enumerate}
\end{dfn}

\begin{example} \label{ex:skewprodcover}
Let $\Lambda$ be a $k$-graph, $G$ a group, $c: \Lambda \to G$ a $1$-cocycle and $G \times_c \Lambda$ the skew product. There is a covering $p : G \times_c \Lambda \to \Lambda$ given by $p(g, \lambda ) = \lambda$ \cite[6.3]{PQR2}. The quotient $G \backslash (G \times_c \Lambda)$ by translation in $G$ is a $k$-graph, and $p$ descends to an isomorphism $\tilde{p} : G \backslash (G \times_c \Lambda) \to \Lambda$.
\end{example}

\begin{thm}[{\cite[Proposition A.19]{BH}, \cite[Theorem 2.7]{PQR2}}] \label{thm:covering universal}
Every connected $k$-graph $\Lambda$ has a universal covering. A connected covering $p :\Sigma \to\Lambda$ is universal if and only if the induced homomorphism $p_*: \pi_1(\Sigma ,v) \to \pi_1 (\Lambda, p(w))$ given by $p_* ( [\gamma] ) = [p(\gamma)]$ is the trivial homomorphism for some, and hence every, $v\in \Sigma^0$.
\end{thm}

\subsection{Simply connected \texorpdfstring{$k$}{k}-graphs}\label{subsec:simplyconnected}

\begin{dfn}
A $k$-graph $\Lambda$ is \emph{simply connected} if $\pi_1(\Lambda, v)$ is trivial for every $v \in \Lambda^0$.
\end{dfn}

\begin{thm}[{\cite[Corollaries 5.5~and~6.5]{PQR2}}]\label{thm:univisspg}
Let $\Lambda$ be a connected $k$-graph.
\begin{enumerate}[(i)]
\item A connected covering $p : \Sigma \to \Lambda$ is universal if and only if $\Sigma$ is simply connected.
\item Given $u \in \Lambda^0$, there exists a cocycle $\eta : \Lambda \to \pi_1 (\Lambda , u)$ for which skew-product covering $p : \pi_1(\Lambda, u) \times_\eta \Lambda  \to \Lambda$ of Example~\ref{ex:skewprodcover} is a universal covering.
\end{enumerate}
\end{thm}

We can characterise simply connected $k$-graphs using either fundamental groupoids or $1$-cocycles.

\begin{lem}\label{lem:characterise simply connected}
Let $\Lambda$ be a connected $k$-graph. Then the following are equivalent:
\begin{enumerate}[(i)]
\item $\Lambda$ is simply connected;
\item $\gamma \mapsto (r(\gamma), s(\gamma))$ is an isomorphism $\Pi(\Lambda) \cong T(\Lambda^0)$; and
\item for every group $G$, every $1$-cocycle $c : \Lambda \to G$ is a coboundary in the sense that there is a function $b : \Lambda^0 \to G$ such that $b(r(\lambda)) c(\lambda) = b(s(\lambda))$.
\end{enumerate}
\end{lem}

\begin{proof}
(i)$\,\Rightarrow\,$(ii). If $\Lambda$ is simply connected, then by definition $\pi_1(\Lambda, v)$ is trivial for all $v$, so the final statement of Remark~\ref{rmk:trans} gives~(ii).
	
(ii)$\,\Rightarrow\,$(iii). Suppose that $\Pi(\Lambda) = T(\Lambda^0)$ and fix a $1$-cocycle $c : \Lambda \to G$. By the universal property of $\Pi(\Lambda)$, there is  a homomorphism $\overline{c} : T(\Lambda^0) \to G$ that extends $c$ (that is, $c = \overline{c}\circ{i}$). Fix $v \in \Lambda^0$. Define $b : \Lambda^0 \to G$ by $b(w) = \overline{c}(v,w)$. For each $\lambda \in \Lambda$
\[
c(\lambda)
= \overline{c}(r(\lambda), s(\lambda))
= \overline{c}((r(\lambda), v) (v, s(\lambda))
= b(r(\lambda))^{-1} b(s(\lambda)),
\]
giving $b(r(\lambda)) c(\lambda) = b(s(\lambda))$.
	
(iii)$\,\Rightarrow\,$(i). Suppose that every $1$-cocycle on $\Lambda$ is a coboundary. Fix $v \in \Lambda^0$. As in Remark~\ref{rmk:trans}, for each $w \in \Lambda^0 \setminus \{v\}$, fix $\gamma_w \in \Pi(\Lambda)^w_v$, and put $\gamma_v = v$. Define $c : \Lambda \to \pi_1(\Lambda, v)$ by $c(\lambda) = \gamma_{r(\lambda)}^{-1} i(\lambda)\gamma_{s(\lambda)}$. Then $c$ is a $1$-cocycle so  there is a map $b : \Lambda^0 \to \pi(\Lambda, v)$ such that $c(\lambda) = b(r(\lambda))^{-1}b(s(\lambda))$ for all $\lambda$.
By the universal property of the fundamental groupoid,  $c$ extends uniquely to a $1$-cocycle $\overline{c} : \Pi(\Lambda) \to \pi_1(\Lambda, v)$ (that is, $c = \overline{c}\circ{i}$).  By uniqueness it follows that for all $\gamma \in \Pi(\Lambda)$, we have
\[
\gamma_{r(\gamma)}^{-1} \gamma\gamma_{s(\gamma)} = \overline{c}(\gamma) = b(r(\gamma))^{-1}b(s(\gamma)).
\]
The first equation implies that the restriction of $\overline{c}$ to $ \pi_1(\Lambda, v)$ is the identity map and by the second equation
the restriction  is trivial.    
Hence, $ \pi_1(\Lambda, v)$ is trivial and so $\Lambda$ is simply connected.
\end{proof}


\subsection{The path groupoid \texorpdfstring{$\Gg_\Lambda$}{} and the \texorpdfstring{$C^*$}{C*}-algebra\texorpdfstring{ $C^*(\Lambda)$}{of a k-graph}}\label{sec:path gpd}

Let $\Lambda$ be a row-finite source-free $k$-graph.  The \emph{infinite path space} $\Lambda^\infty$ of $\Lambda$ is the space of $k$-graph morphisms $x : \Omega_k \to \Lambda$ under the locally compact Hausdorff topology with basic compact open sets $Z(\lambda) := \{ x \in  \Lambda^\infty: \lambda = x(0, d(\lambda)) \}$, indexed by $\lambda \in \Lambda$. For $p \in \NN^k$, the shift map $\sigma^p: \Lambda^\infty \to \Lambda^\infty$ is defined by $\sigma^px(m, n) = x(m+p, n+p)$ for $x \in \Lambda^\infty$ and $(m, n) \in \Omega_k$; and $p \mapsto \sigma^p$ is an action of $\NN^k$ by local homeomorphisms.

Elements $x, y \in \Lambda^\infty$ are \emph{shift equivalent}, written $x \simeq y$, if $\sigma^px = \sigma^qy$ for some $p, q \in \NN^k$. The path groupoid $\Gg_\Lambda$ is the Deaconu--Renault groupoid of the action $p \mapsto \sigma^p$:

\begin{dfn}[{\cite[Definition~2.7]{KP2}}]\label{dfn:pathgpd}
	The path groupoid is
	\[
	\Gg_\Lambda := \{ (x, n, y) \in \Lambda^\infty \times \ZZ \times \Lambda^\infty : \sigma^\ell x = \sigma^my, n = \ell -m \},
	\]
	with unit space $\Gg_\Lambda^0 = \{(x, 0, x) : x \in \Lambda^\infty\}$ identified with $\Lambda^\infty$,
	with structure maps
	\[
	r(x, n, y) = x, \quad s(x, n, y) = y, \quad (x, n, y)(y, \ell, z) = (x, n+\ell, z), \quad (x, n, y)^{-1} =  (y, -n, x),
	\]
    and under the topology with basic open sets $Z(\mu,\nu) = \{(\mu x, d(\mu) - d(\nu), \nu x) : x \in Z(s(\mu))\}$ indexed by pairs $(\mu,\nu) \in \Lambda \times \Lambda$ such that $s(\mu) = s(\nu)$.
\end{dfn}

The $C^*$-algebra of $\Lambda$ is defined via generators and relations.

\begin{dfn}[{\cite[Definitions~1.5]{KP2}}]
	A family of partial isometries $\{ s_\lambda :  \lambda \in \Lambda \}$ is a \emph{Cuntz--Krieger $\Lambda$-family} if
	\begin{itemize}
		\item[(CK1)] $\{ s_v : v \in  \Lambda^0\}$ is a collection of mutually orthogonal projections;
		\item[(CK2)] $s_{\lambda\mu} = s_{\lambda}s_\mu$ for all $\lambda, \mu \in \Lambda$ such that $s(\lambda) = r(\mu)$;
		\item[(CK3)] $s_{\lambda}^*s_{\lambda}  = s_{s(\lambda)}$ for all $\lambda \in \Lambda$; and
		\item[(CK4)] for all $v \in \Lambda^0$ and $n \in \NN^k$ we have $s_v = \sum_{\lambda \in v\Lambda^n}  s_{\lambda}s_{\lambda}^*$.
	\end{itemize}
	We write $C^*(\Lambda)$ for the universal  $C^*$-algebra generated by a Cuntz--Krieger $\Lambda$-family $\{ s_\lambda :  \lambda \in \Lambda \}$.
\end{dfn}

The groupoid $\Gg_\Lambda$ is \'etale \cite[Proposition~2.8]{KP2}, and $C^*(\Lambda) \cong C^*(\Gg_\Lambda)$ \cite[Corollary~3.5(i)]{KP2}.

\section{Embedding results for higher-rank graphs} \label{sec:bootstrap}

In this section we develop tools for determining when a $k$-graph $\Lambda$ embeds in $\Pi(\Lambda)$, and describe classes of examples that do embed; we also present three examples---one from \cite{PQR1}, one due to Ben Steinberg, and one that is new---that do not embed.

\subsection{Non-embeddings}\label{subsec:non-embeddings}

Even a fairly elementary monoidal $2$-graph $\Lambda$ need not embed in $\Pi(\Lambda)$:

\begin{example}[{\cite[Example~7.1]{PQR1}}]\label{ex:pqrex1}
Let $\Lambda$ be the 1-vertex $2$-graph with $\Lambda^{\varepsilon_1} = \{ d, e \}$ and $\Lambda^{\varepsilon_2} = \{ a,b,c \}$ such that
\begin{equation} \label{eq:pqrex}
da = ad,\quad   db = be,\quad  dc = ae,\quad  ea = cd,\quad   eb  = ce,\quad  ec = bd.
\end{equation}
Using the first four relations from \eqref{eq:pqrex} and that the map $i : \Lambda \to \Pi (\Lambda)$ is a morphism we obtain
	\begin{align*}
		i(a) =  i(d)  i(a) i(d)^{-1} &= i(d)i(e)^{-1}i(c) = i(d)i(b)i(e)^{-1} = i(b),
	\end{align*}
so $i(a)=i(b)$ in $\Pi ( \Lambda)$. The fifth equation in~\eqref{eq:pqrex} gives $i(d)=i(e)$, so equations two and five give $i(b)=i(c)$. Hence $i(a)=i(b) = i(c)$ and $i(d) = i(e)$. The degree map descends to an isomorphism $\tilde{d} : \Pi ( \Lambda ) \to \ZZ^2$; so the universal cover of $\Lambda$ is isomorphic to $\ZZ^2 \times_d \Lambda$.
\end{example}

The next example, shown to us by Ben Steinberg, who attributes the idea to Ma{l}'cev (see also \cite[Example 11.13]{LV}), is a monoidal $2$-graph that does not embed even though its edge-set does.

\begin{example}[Steinberg, private communication]  \label{ex:Steinberg}
Let $\Lambda$ be the unique $1$-vertex $2$-graph with $\Lambda^{\varepsilon_1} = \{ e_1, e_2, e_3, e_4 \}$ and $\Lambda^{\varepsilon_2} = \{ f_1, f_2, f_3, f_4 \}$, and such that
\begin{equation}\label{eq:steinbergfactorisationrules}
	e_af_b :=
\begin{cases}
		f_be_a & \text{if  } (a, b) = (1, 4), (4, 1) ; \\
		f_ae_b & \text{otherwise.}
\end{cases}
\end{equation}
Since $i : \Lambda \to \Pi(\Lambda)$ is a functor, $i(e_a)^{-1}i(f_a) = i(f_b)i(e_b)^{-1}$ for $(a, b) = (1,2), (3,2), (3,4)$, so
\[
i( e_1 )^{-1} i ( f_1) = i( f_2 ) i ( e_2 )^{-1} = i( e_3 )^{-1} i ( f_3) = i ( f_4 ) i (e_4 )^{-1},
\]
and then rearranging the outer terms gives
\begin{equation}\label{eq:extra rel}
i(f_1e_4) = i(f_1) i(e_4) =  i(f_4) i (e_1) = i(f_4e_1).
\end{equation}
Uniqueness of factorisations in $\Lambda$ shows that $f_1e_4 \not= f_4e_1$, so $i$ is not injective.

We show that $i$ is injective on $\Lambda^{\varepsilon_1} \cup \Lambda^{\varepsilon_2}$. For this, define $c : \Lambda^{\varepsilon_1} \cup \Lambda^{\varepsilon_2} \to \ZZ$ by $c(e_j) = c(f_j) = j$ for $j=1, \ldots , 4$. Since $c$ respects the equations~\eqref{eq:steinbergfactorisationrules}, it extends to a functor $c : \Lambda \to \ZZ$. By Definition~\ref{dfn:fundgpd} there is a functor $\tilde{c} : \Pi(\Lambda)
\to \ZZ$ such that $\tilde{c} \circ i = c$. In particular, $\tilde{c}(i(e_j)) = j = \tilde{c}(i(f_j))$ for all $j$. Hence $(\tilde{c} \times \tilde{d}) \circ i : (\Lambda^{\varepsilon_1} \cup \Lambda^{\varepsilon_2}) \to \ZZ^2$ is injective. Thus $i$ is injective on $\Lambda^{\varepsilon_1} \cup \Lambda^{\varepsilon_2}$.
\end{example}

\begin{example}\label{ex:Davids3graph}
By \cite[Theorems 4.4~and~4.5]{HRSW}, there is a unique $3$-graph $\Gamma$ with the skeleton and factorisation rules below (there are no 3-coloured paths, so the associativity condition is vacuous).
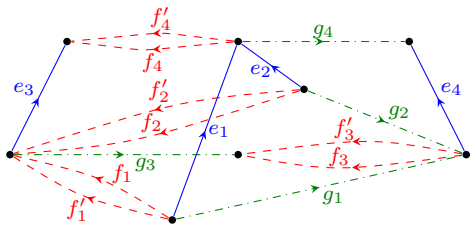
\[
\begin{tikzpicture}[scale=1.5, >=stealth, decoration={markings, mark=at position 0.5 with {\arrow{>}}}]
    \node[circle, inner sep=1pt, fill=black] (v0) at (0,0,1.5) {};
    \node[circle, inner sep=1pt, fill=black] (v1) at (0,0,-1.5) {};
    \node[circle, inner sep=1pt, fill=black] (v2) at (0,0,0) {};
    \node[circle, inner sep=1pt, fill=black] (v3) at (2,0,0) {};
    \node[circle, inner sep=1pt, fill=black] (v4) at (-2,0,0) {};
    \node[circle, inner sep=1pt, fill=black] (v5) at (1.5,1,0) {};
    \node[circle, inner sep=1pt, fill=black] (v6) at (-1.5,1,0) {};
    \node[circle, inner sep=1pt, fill=black] (v7) at (0,1,0) {};
    \draw[red, dashed, postaction=decorate] (v1) .. controls +(-1,0,0.25) and +(1,0,-1.25) .. (v4) node[red, pos=0.5, anchor=south east, inner sep=1pt] {\tiny$f'_2$};
    \draw[red, dashed, postaction=decorate] (v1) .. controls +(-1,0,1.25) and +(1,0,-0.25) .. (v4) node[red, pos=0.45, anchor=south, inner sep=1.2pt] {\tiny$f_2$};
    \draw[green!50!black, dash dot, postaction=decorate] (v1)--(v3) node[green!50!black, pos=0.5, anchor=south west, inner sep=1pt] {\tiny$g_2$};
    \draw[red, dashed, postaction=decorate] (v0) .. controls +(-1,0,-0.5) and +(1,0,1) .. (v4) node[red, pos=0.45, anchor=north east, inner sep=0pt] {\tiny$f'_1$};
    \draw[red, dashed, postaction=decorate] (v0) .. controls +(-1,0,-1) and +(1,0,0.5) .. (v4) node[red, pos=0.4, anchor=south west, inner sep=0pt] {\tiny$f_1$};
    \draw[green!50!black, dash dot, postaction=decorate] (v0)--(v3) node[green!50!black, pos=0.5, anchor=north west, inner sep=1pt] {\tiny$g_1$};
    \draw[red, dashed, bend right, postaction=decorate] (v3) .. controls +(-1,0,-0.4) and +(1,0,-0.4) .. (v2) node[red, pos=0.62, anchor=south, inner sep=0.5pt] {\tiny$f'_3$};
    \draw[red, dashed, bend left,  postaction=decorate] (v3) .. controls +(-1,0,0.4) and +(1,0,0.4) .. (v2) node[red, pos=0.5, anchor=south, inner sep=0.5pt] {\tiny$f_3$};
    \draw[green!50!black, dash dot, postaction=decorate] (v4)--(v2) node[green!50!black, pos=0.6, anchor=north, inner sep=1pt] {\tiny$g_3$};
    \draw[red, dashed, bend right, postaction=decorate] (v7) .. controls +(-1,0,-0.25) and +(1,0,-0.25) .. (v6) node[red, pos=0.5, anchor=south, inner sep=1.2pt] {\tiny$f'_4$};
    \draw[red, dashed, bend left,  postaction=decorate] (v7) .. controls +(-1,0,0.25) and +(1,0,0.25) .. (v6) node[red, pos=0.4, anchor=north, inner sep=1.2pt] {\tiny$f_4$};
    \draw[green!50!black, dash dot, postaction=decorate] (v7)--(v5) node[green!50!black, pos=0.5, anchor=south, inner sep=1pt] {\tiny$g_4$};
    \draw[blue, postaction=decorate] (v1)--(v7) node[blue, pos=0.5, anchor=north east, inner sep=-0.5pt] {\tiny$e_2$};
    \draw[blue, postaction=decorate] (v0)--(v7) node[blue, pos=0.5, anchor=west, inner sep=1pt] {\tiny$e_1$};
    \draw[blue, postaction=decorate] (v3)--(v5) node[blue, pos=0.5, anchor=south west, inner sep=1pt] {\tiny$e_4$};
    \draw[blue, postaction=decorate] (v4)--(v6) node[blue, pos=0.5, anchor=south east, inner sep=1pt] {\tiny$e_3$};
    \node at (4.8,0.25,0) {\parbox[t]{7cm}{\small${e_3}{f_i} = {f_4}{e_i}$ \quad ${e_3}{f'_i} = {f'_4}{e_i}$\quad($i=1,2$)\\
                                           ${e_4}{g_i} = {g_4}{e_i}$\phantom{ \quad ${e_1}{f'_i} = {f'_4}{e_i}$}\quad($i=1,2$)\\
                                           ${f_3}{g_2} = {g_3}{f_2}$ \quad ${f_3'}{g_2} = {g_3}{f'_2}$\\
                                           ${f'_3}{g_1} = {g_3}{f_1}$ \quad ${f_3}{g_1} = {g_3}{f_1'}$
                     }};
\end{tikzpicture}
\]
By Proposition~\ref{lprop:ess-emb}(ii) below, the 2-coloured sub-$2$-graphs of $\Gamma$ are all embeddable: the $\ZZ$-valued cocycle on the blue-red graph carrying $\{f_i : i \le 4\}$ to $1$ and all other edges to $0$ is essential (see Definition~\ref{def:essential}); the $\ZZ$-valued cocycle on the red-green graph carrying $\{f_1\}\cup\{f'_i :  i \ge 2\}$ to $1$ and all other edges to $0$ is essential; and the trivial cocycle on the blue-green graph is essential.

However $\Gamma$ does not embed in $\Pi(\Gamma)$: writing $[x]$ for $i(x) \in \Pi(\Lambda)$, we calculate:
\begin{align*}
[f_1'] [g_1]^{-1}
    &= [g_3]^{-1}[f_3]
    = [f_2] [g_2]^{-1}
    = [f_2] [e_2]^{-1} [e_2] [g_2]^{-1}
    = [e_3]^{-1}[f_4] [g_4]^{-1}[e_4] \\
    &= [e_3]^{-1}[f_4][e_1][e_1]^{-1} [g_4]^{-1}[e_4]
    = [e_3]^{-1}[e_3][f_1][g_1]^{-1} [e_4]^{-1}[e_4]
    = [f_1] [g_1]^{-1}.
\end{align*}
So cancellation gives $[f_1'] = [f_1]$. We then have $[f'_3] = [g_3][f_1][g_1]^{-1} = [g_3][f'_1][g_1]^{-1} = [f_3]$ and $[f'_4] = [e_3][f'_1][e_1]^{-1} = [e_3][f_1][e_1]^{-1} = [f_4]$, and then also $[f'_2] = [e_3]^{-1} [f'_4][e_2] = [e_3]^{-1} [f_4][e_2] = [f_2]$.
\end{example}

Motivated by these examples, we seek conditions under which $i : \Lambda \to \Pi(\Lambda)$ \emph{is} injective.

\subsection{Embedding singly connected higher-rank graphs}\label{subsec:singlyconnected}

\begin{dfn}
A  $k$-graph $\Lambda$ is \emph{singly connected} if there is at most one path between any two  vertices; that is, for all $u,v \in \Lambda^0$ we have $| u {\Lambda}v | \le 1$.
\end{dfn}

Singly connected $k$-graphs need not be connected. The vertex set of a singly-connected $k$-graph is partially ordered by the relation $\le$ given by $u \le v$ if and only if $u \Lambda v \neq \emptyset$.

\begin{example} \label{ex:essential-cocycle}
Write $\{ t_i : i=1 , \ldots n\}$ for the generators of the free group $\FF_n$. Let $c :B_n \to \FF_n$ be the $1$-cocycle such that $c(f_i) = t_i$ for all $i$. Then $\FF_n \times_c B_n$ is singly connected.
\end{example}

There is a relationship between singly connected $k$-graphs, and the simply connected $k$-graphs of Section~\ref{subsec:simplyconnected}, though neither condition implies the other.

\begin{prop}\label{prop:1cnc}
Let $\Lambda$ be a connected $k$-graph and suppose that $i : \Lambda \to \Pi(\Lambda)$ is injective.
If $\Lambda$ is simply connected, then it is singly connected.
\end{prop}

\begin{proof}
Suppose that $\Lambda$ is not singly connected.
Then there exist distinct elements $\lambda, \mu \in \Lambda$ such that $s(\lambda) = s(\mu)$ and  $r(\lambda) = r(\mu)$.
Since $i : \Lambda \to \Pi(\Lambda)$ is injective,  $i(\lambda) \ne i(\mu)$ and thus  $i(\lambda)^{-1}i(\mu) \in \pi_1 (\Lambda, s(\lambda)) \backslash \{ s(\lambda)\}$.  Hence, $\Lambda$ is not simply connected.
\end{proof}

The reverse implication fails, as the following example illustrates.

\begin{example}
Let $E$ be the directed graph with $E^0=\{u, v, w, x\}$ and $E^1= \{ e, f, g, h \}$ such that $s(e) = u = s(f)$,  $s(g) = w = s(h)$, $r(e) = v = r(h)$ and $r(f) = x = r(g)$. Then the $1$-graph $E^*$ is a singly-connected $1$-graph that is not simply connected since $\pi_1(E^*, u) \cong \ZZ$.   Adding tails at both $u$ and $w$ as in \cite[Lemma~1.2]{BPRS} yields a source-free $1$-graph with the same property.
\end{example}

We will use the next theorem, which exploits the universal property of the fundamental groupoid from Section~\ref{sec:fun}, to show that singly connected $k$-graphs embed in their fundamental groupoids.
\begin{thm} \label{prop:prec}
	Let $\Lambda$ be a $k$-graph and let $\Gg$ be a groupoid. If there is an injective functor $F : \Lambda \to \Gg$,
	then $i : \Lambda \to \Pi(\Lambda)$ is injective.
\end{thm}

\begin{proof}
	The universal property of the fundamental groupoid, yields a homomorphism $\tilde{F} : \Pi(\Lambda) \to  \Gg$
	such that $F = \tilde{F} \circ i$.
	Hence if $F$ is injective, then $i$ is injective.
\end{proof}

\begin{prop}\label{prop:1-cnc-emb}
	Let $\Lambda$ be a connected $k$-graph.  Then
	\begin{enumerate}[(i)]
		\item  The canonical map  $\iota: \Lambda \to T(\Lambda^0)$  is injective if and only if $\Lambda$ is singly connected.
		\item If $\Lambda$ is singly connected, then $i : \Lambda \to \Pi(\Lambda)$ is injective.
	\end{enumerate}
\end{prop}

\begin{proof}
	The first assertion follows by definition and the second follows from Theorem~\ref{prop:prec}.
\end{proof}

Theorem~\ref{prop:prec} also allows us to deduce embeddability from the existence of a suitable $1$-cocycle.

\begin{dfn}\label{def:essential}
Let  $\Lambda$ be a $k$-graph, $G$ a countable group, and $c : \Lambda \to G$ a 1-cocycle. We say that $c$ is \emph{essential} if the restriction of $c$ to $u{\Lambda}v$ is injective for all $u, v \in \Lambda^0$.
\end{dfn}

\begin{example}
The $1$-cocycle $c : B_n \to \FF_n$ described in Example~\ref{ex:essential-cocycle} is essential.
\end{example}

\begin{prop}\label{lprop:ess-emb}
Let  $\Lambda$ be a $k$-graph. Then the following are equivalent:
\begin{enumerate}[(i)]
\item the canonical cocycle $\kappa : \Lambda \to \pi_1(\Lambda)$ is essential;
\item $\Lambda$ admits an essential cocycle  $c : \Lambda \to G$ to a group $G$; and
\item  $i : \Lambda \to \Pi(\Lambda)$ is injective.
\end{enumerate}
For any essential cocycle $c : \Lambda \to G$ as in~(ii), $G \times_c \Lambda$ is singly connected.
\end{prop}

\begin{proof}
The implication (i)~$\implies$~(ii) is obvious.

For (ii)~$\implies$~(iii), suppose that $c : \Lambda \to G$ is an essential cocycle into a group.
Note that $G  \times T(\Lambda^0)$ is a groupoid.
Define $j : \Lambda \to G \times T(\Lambda^0)$  by $j(\lambda) := (c(\lambda), (r(\lambda), s(\lambda)))$; then $j$ is a functor.
Since $c$ is essential, $j$ is injective, so $i : \Lambda \to \Pi(\Lambda)$ is injective by Theorem~\ref{prop:prec}.

For (iii)~$\implies$~(i), suppose that  $i : \Lambda \to \Pi(\Lambda)$ is injective.
Fix $u,v \in \Lambda^0$ and $\mu \not= \nu$ in $v\Lambda u$. For each $w \in [v] = s(\Pi(\Lambda)^v)$, choose $\gamma_w \in \Pi(\Lambda)^v_w$, with $\gamma_v = \{v\}$. Define $c : \Lambda \to \pi_1(\Lambda, v)$ by
\[
c(\lambda) = \begin{cases}
\gamma_{r(\lambda)}i(\lambda)\gamma_{s(\lambda)}^{-1} &\text{ if $r(\lambda) \in [v]$}\\
e &\text{ otherwise.}
\end{cases}
\]
Then $c$ is a cocycle. We have $c(\mu)c(\nu)^{-1} = \gamma_v i(\mu) \gamma_u^{-1} \gamma_u i(\nu) \gamma_v^{-1} = \gamma_vi(\mu)i(\nu)^{-1}\gamma_v^{-1}$.
Since $i$ is injective, we obtain $c(\mu) \not= c(\nu)$; and then the universal property of $\kappa : \Lambda \to \Pi(\Lambda)$
in Proposition~2.13 implies that $\kappa(\mu) \not= \kappa(\nu)$.  Hence, $\kappa : \Lambda \to \pi_1(\Lambda)$ is essential.

For the final statement, suppose that $c : \Lambda \to G$ is essential, and that $r(g, \lambda) = r(h,\mu)$ and $s(g,\lambda) = s(h,\mu)$ in $G \times_c \Lambda$.
Then
\begin{align*}
(g, r(\lambda)) &= r(g,\lambda) = r(h, \mu) = (h, r(\mu)), \quad \text{and} \\
(c(\lambda) g, s(\lambda)) &= s(g,\lambda) = s(h,\mu) = (c(\mu) h, s(\mu)).
\end{align*}
So $r(\lambda) = r(\mu)$, $s(\lambda) = s(\mu)$,
$g = h$ and $c(\lambda)g = c(\mu)h$. These last two equalities give $c(\lambda) = c(\mu)$.
Thus $j(\lambda) = j(\mu)$ and hence $\lambda = \mu$.
Therefore $(g, \lambda) = (h,\mu)$ and so  $G \times_c \Lambda$ is singly connected.
\end{proof}

\subsection{More-general embedding results}\label{subsec:generalembeddingresults}

In this section we investigate embeddability of $k$-graphs that are not singly connected. We start with one of the most useful results in our toolkit, which relies on the universal property of the fundamental groupoid given in Definition~\ref{dfn:fundgpd}.

\begin{prop}[Lifting embeddability] \label{prop:lift}
Let $\Lambda , \Sigma$ be connected $k$-graphs and let $p : \Sigma \to \Lambda$ be a covering. Then $i_\Lambda : \Lambda \to \Pi ( \Lambda )$ is injective if and only if $i_\Sigma : \Sigma \to \Pi(\Sigma)$ is  injective.
\end{prop}
\begin{proof}
Suppose that $i_\Lambda : \Lambda \to \Pi(\Lambda )$ is injective and that $\sigma, \sigma'  \in \Sigma$ satisfy $ i_\Sigma(\sigma) = i_\Sigma( \sigma' )$.
In particular, $s(\sigma) = s(\sigma')$; let $u := s(\sigma)$.
By universality of $\Pi(\Sigma)$, there is a unique groupoid morphism $\tilde{p} : \Pi(\Sigma) \to \Pi(\Lambda)$ such that
$\tilde{p} \circ  i_\Sigma = i_\Lambda \circ p$.   Hence
\[
 i_\Lambda(p(\sigma)) = \tilde{p} ( i_\Sigma( \sigma)) =  \tilde{p} ( i_\Sigma( \sigma')) =  i_\Lambda(p(\sigma')).
\]
Injectivity of $i_\Lambda$ forces $p(\sigma) = p(\sigma')$.  Since $p$ is a covering, it is injective on $s^{-1}(u)$.
So $\sigma = \sigma'$ and hence $i_\Sigma : \Sigma \to \Pi(\Sigma)$, is injective.

For the reverse implication suppose that  $i_\Lambda : \Lambda \to \Pi ( \Lambda )$ is  not injective. Then there are  distinct $\lambda, \lambda' \in \Lambda$  such that $i_\Lambda( \lambda ) = i_\Lambda( \lambda' )$.
We may assume without loss of generality that $\Sigma$ is the universal covering of $\Lambda$ so that $\Sigma$ is simply connected. Since $\Sigma$ is connected, $r \times s : \Pi ( \Sigma ) \to \Sigma^0 \times \Sigma^0$ is an isomorphism, so $\Pi ( \Sigma) \cong T( \Sigma^0 )$. By Theorem~\ref{thm:univisspg}(ii), given $u \in \Lambda^0$,
there is a cocycle $\eta : \Lambda \to \pi_1(\Lambda, u)$ such that $\Sigma \cong \pi_1(\Lambda, u) \times_\eta \Lambda$
and $p$ is given by projection onto the second factor. It follows that $s( \lambda ) = s( \lambda' )$, $r( \lambda ) = r( \lambda' )$ and  $\eta( \lambda ) = \eta( \lambda' )$
(since $\eta$ factors through $i_\Lambda$ and $i_\Lambda( \lambda ) = i_\Lambda( \lambda' )$).
Identifying $\Sigma$ with the skew-product as above, set $\sigma = (1, \lambda), \sigma' =  (1, \lambda') \in \Sigma$, so $\sigma, \sigma'$ are distinct. We have
\[
s(\sigma) = s(1, \lambda) = (\eta(\lambda), s(\lambda)) = (\eta(\lambda'), s(\lambda')) = s(1, \lambda') = s(\sigma')
\]
and  similarly $r(\sigma) = r(\sigma')$. So $(r \times s)(i_\Sigma(\sigma)) = (r \times s)(i_\Sigma(\sigma'))$. Since $r \times s$ is injective on $i_\Sigma$, we deduce that $i_\Sigma : \Sigma \to \Pi(\Sigma )$ is  not injective.
\end{proof}

Our later results say that embeddability is preserved by various constructions of new $k$-graphs from old ones. So we need to know that some basic classes of $k$-graphs, like $1$-graphs, embed.

\begin{thm} \label{thm:1gcase}
Let $\Lambda$ be a $1$-graph. Then $i_\Lambda : \Lambda \to \Pi ( \Lambda )$  is injective.
\end{thm}
\begin{proof}  
Write $\Lambda = \bigsqcup_{i=1}^n \Lambda_i $ as a disjoint union of connected graphs. For $i=1 , \ldots , n$ let $\Sigma_i$ be the universal cover of $\Lambda_i$.
Since $\Sigma = \bigsqcup_{i=1}^n \Sigma_i$ is (the path category of) a disjoint union of trees, there is at most one undirected path connecting any two distinct vertices.
It follows that $\Sigma$  is singly connected and therefore embeddable by Proposition~\ref{prop:1-cnc-emb}(ii).
Hence, $i_\Lambda : \Lambda \to \Pi ( \Lambda )$  is injective.
\end{proof}

\begin{cor}\label{cor:specialcases}
Let $\Lambda$ be a $k$-graph, and suppose that $i_\Lambda : \Lambda \to \Pi(\Lambda)$ is injective.
\begin{enumerate}[(i)]
\item\label{it:affinepullbackcase} Let $f : \NN^\ell \to \NN^k$ be an affine map. Then $i : f^* ( \Lambda) \to \Pi(f^* ( \Lambda ) )$ is injective.
\item\label{it:cartprodcase} If $\Gamma$ is an $\ell$-graph and $i_\Gamma : \Gamma \to \Pi(\Gamma)$ is injective, then $i_{\Lambda
\times \Gamma} : \Lambda \times \Gamma \to \Pi(\Lambda \times \Gamma)$ is injective.
\item\label{it:skewprodcase} If $c : \Lambda \to G$ is a $1$-cocycle into a group, then $i_{G\times \Lambda} : G \times_c \Lambda \to \Pi ( G \times_c \Lambda )$ is injective.
\item\label{it:crosscase} If $\alpha : \NN^\ell \to \operatorname{Aut} (\Lambda)$ is an action, then there is an action $\tilde{\alpha} : \ZZ^\ell \to \Aut(\Pi(\Lambda))$ such that $\tilde{\alpha_n} \circ i_\Lambda = i_{\Lambda} \circ \alpha_n$ for $n \in \NN^\ell$. Both
\[
i_\Lambda \times i_{\NN^\ell} : \Lambda \times_\alpha  \NN^\ell \to \Pi ( \Lambda) \times_\alpha \ZZ^\ell \quad\text{and}\quad
i_{\Lambda \times_\alpha \NN^\ell} : \Lambda \times_\alpha \NN^\ell \to \Pi(\Lambda \times_\alpha \NN^\ell)
\]
are injective. Moreover, $i_\Lambda \times i_{\NN^\ell}$ induces an isomorphism $\Pi ( \Lambda \times_\alpha \NN^\ell ) \cong \Pi ( \Lambda) \times_\alpha \ZZ^\ell$.
\end{enumerate}
\end{cor}
\begin{proof}
(\ref{it:affinepullbackcase}) Define $i \times \id : f^*( \Lambda ) \to \Pi(\Lambda) \times \ZZ^\ell$ by $(i \times \id)(\lambda, n) = (i(\lambda), n)$. Then $i \times \id$ is an injective functor into a groupoid, so the result follows from Theorem~\ref{prop:prec}.

(\ref{it:cartprodcase}) The map $i_\Lambda \times i_\Gamma : \Lambda \times \Gamma \to \Pi(\Lambda) \times \Pi(\Gamma)$ is an injective functor into a groupoid, so the result
follows from Theorem~\ref{prop:prec}.

(\ref{it:skewprodcase}) By universality of $\Pi(\Lambda)$, there is a cocycle $\tilde{c} : \Pi(\Lambda) \to G$ such that $\tilde{c} \circ i_\Lambda = c$. The skew-product groupoid $G \times_{\tilde{c}} \Pi(\Lambda)$ is equal as a set to $G \times \Pi(\Lambda)$, and $\id_G \times i_\Lambda : G \times_c \Lambda \to G \times_{\tilde{c}} \Pi(\Lambda)$ is a functor. Since $i_\Lambda$ is injective, so is $\id_G \times i_\Lambda$, so the result follows from Theorem~\ref{prop:prec}.

(\ref{it:crosscase}) Since the action $\alpha$ of $\NN^\ell$ on $\Lambda$ is determined by $\ell$ commuting automorphisms, it extends to an action (also called $\alpha$) of $\ZZ^\ell$ on $\Lambda$. By functoriality this extends to an action $\tilde{\alpha} : \ZZ^\ell \to \Aut(\Pi(\Lambda))$ such that $\tilde{\alpha_n} \circ i_\Lambda = i_{\Lambda} \circ \alpha_n$ for $n \in \NN^\ell$. It is routine to check that $i_\Lambda \times i_{\NN^\ell}$ is a functor; it is injective because $i_\Lambda$ and $i_{\NN^\ell}$ are injective. So Theorem~\ref{prop:prec} implies that $i_{\Lambda \times_\alpha \NN^\ell}$ is injective.

To see that $i_\Lambda \times i_{\NN^\ell}$ induces an isomorphism $\Pi ( \Lambda \times_\alpha \NN^\ell ) \cong \Pi ( \Lambda) \times_\alpha \ZZ^\ell$, note that the universal property of $\Pi ( \Lambda \times_\alpha \NN^\ell )$ implies that $i_\Lambda \times i_{\NN^\ell}$ induces a homomorphism $\tilde{i} : \Pi ( \Lambda \times_\alpha \NN^\ell ) \to \Pi ( \Lambda) \times_\alpha \ZZ^\ell$ such that $\tilde{i} \circ i_{\Lambda \times_\alpha \NN^\ell} = i_\Lambda \times i_{\NN^\ell}$. We construct an inverse.
The restriction $c_1 := i_{\Lambda \times_\alpha \NN^\ell}|_{\Lambda \times \{0\}} : \Lambda \to \Pi(\Lambda \times_\alpha \NN^\ell)$ is a functor, as is $c_2 := i_{\Lambda \times_\alpha \NN^\ell}|_{\Lambda^0 \times \NN^\ell}$. The universal property of $\Pi(\Lambda)$ implies that $i_{\Lambda \times_\alpha \NN^\ell}|_{\Lambda \times \{0\}}$ induces a homomorphism $\tilde{c}_1 : \Pi(\Lambda) \to \Pi ( \Lambda \times_\alpha \NN^\ell )$; and $i_{\Lambda \times_\alpha \NN^\ell}|_{\Lambda^0 \times \NN^l}$ extends to a homomorphism $\tilde{c}_2 : \Lambda^0 \times \ZZ^\ell \to \Pi ( \Lambda \times_\alpha \NN^\ell )$. Routine calculations show that $\tilde{c}_1 \times \tilde{c}_2 : \Pi(\Lambda) \times_{\tilde{\alpha}} \ZZ^\ell \to \Pi(\Lambda \times_\alpha \NN^\ell)$ is a homomorphism inverse to $\tilde{i}$.
\end{proof}

\begin{rmk}\label{rmk:1graphcp embeds}
Combining Theorem~\ref{thm:1gcase} and Corollary~\ref{cor:specialcases}(\ref{it:crosscase}), we see that crossed-product graphs of $1$-graphs always embed in their fundamental groupoids.
\end{rmk}

\begin{examples} \label{ex:thetaswap}
We present two examples of Corollary~\ref{cor:specialcases}(\ref{it:affinepullbackcase}).
\begin{enumerate}[(i)]
\item Define $f : \NN^2 \to \NN$ by  $f(a,b)=a+b$, and let $\Lambda = f^* (B_n)$.  By Corollary~\ref{cor:specialcases}(\ref{it:affinepullbackcase}) $f^* (\Lambda)$ embeds in its fundamental group since the $1$-graph $B_n$ does by Theorem~\ref{thm:1gcase}. Indeed, for $\theta : [n] \times [n] \to [n] \times [n]$ given by $\theta (i,j)=(i, j)$, we have $\Lambda \cong \mathbb{F}_\theta^+$.
\item Let $\Lambda$ be a $2$-graph and define $f : \NN^2 \to \NN^2$ by $f(a,b)=(a,b)+{\bf 1}$. Then $f^*(\Lambda)$ is the dual graph ${\bf 1} \Lambda$ described in \cite[Definition 3.1]{APS}. So for the $2$-graph $\Lambda = \FF_\theta^+$ from~(i) above, $f^* ( \FF_\theta^+ )$ embeds in its fundamental group by Corollary~\ref{cor:specialcases}(\ref{it:affinepullbackcase}).
\end{enumerate}
\end{examples}

\begin{cor}[Action graphs] \label{thm:genaction}
Let $\Lambda$ be a $k$-graph. Let $B_n$ be the $1$-graph described in Examples~\ref{ex:bnex}(i). Let $\mu \mapsto \alpha_\mu$ be a functor from $B_n$ to $\Aut(\Lambda)$. Let $\Gamma = B_n \times \Lambda$; define $d : \Gamma \to \NN^{k+1}$ by $d(\mu, \lambda) = (|\mu|, d(\lambda))$; define $r, s : \Gamma \to \Gamma^0$ by $r(\mu, \lambda) = (u, \alpha_\mu(r(\lambda )))$ and $s(\mu, \lambda) = (u, s(\lambda))$; and define composition in $\Gamma$ by
\begin{equation} \label{eq:gammadef}
(\mu,\alpha_\nu(\lambda))(\nu,\xi) = (\mu\nu, \lambda\xi) .
\end{equation}
Then $(\Gamma,d)$ is a $(k+1)$-graph. If $i_\Lambda : \Lambda \to \Pi ( \Lambda )$ is injective, then $i_\Gamma : \Gamma \to \Pi ( \Gamma )$ is injective.
\end{cor}

\begin{proof}
It is routine to check that~\eqref{eq:gammadef} determines an associative composition. The map $d$ is clearly a functor, and if $d(\mu,\lambda) = (a+b, m+n)$, then factorising $\mu = \mu_a \mu_b$ and $\lambda = \lambda_m \lambda_n$ with the appropriate degrees, the factorisation $(\mu,\lambda) = (\mu_a, \alpha_{\mu_b}(\lambda_m))(\mu_b,\lambda_n)$ is the unique factorisation of $(\mu,\lambda)$ into morphisms of degree $(a,m)$ and $(b,n)$. So $\Gamma$ is a $(k+1)$-graph.

Universality of $\Pi ( \Lambda)$ implies that each $\alpha_\mu$ extends to an automorphism of $\Pi ( \Lambda )$. So $\alpha$ extends to an action of $\Pi(B_n) \cong \mathbb{F}_n$ on $\Pi ( \Lambda )$, with semidirect product groupoid $\Pi(B_n) \ltimes_{\tilde{\alpha}} \Pi ( \Lambda )$. Then $i_{B_n} \times i_\Lambda$ is an embedding of $\Gamma$ in $\Pi(B_n) \ltimes_{\tilde{\alpha}} \Pi ( \Lambda )$. The result now follows by Theorem~\ref{prop:prec}.
\end{proof}

\begin{examples} \label{ex:bnex3}
\begin{enumerate}[(i)]
\item Fix $m,n  \ge 2$. Let $\Lambda = B_m$ be the  $1$-graph described in Examples~\ref{ex:bnex}(i). For each $f \in B_n^1$, let $\alpha_f$ be a permutation of $B_m^1$, and extend this to a $1$-cocycle $B_n \to \Aut(\Lambda)$ in the only possible way. By Theorem~\ref{thm:genaction} this data gives rise to a $2$-graph $\Gamma$ that embeds in its fundamental group. Define $\theta : [n] \times [m] \to [m] \times [n]$ by $\theta(i,j) = (j' , i)$ if and only if $\alpha_{f_i}(f_j) = f_{j'}$. Then $\Gamma$ is isomorphic to the $2$-graph $\mathbb{F}_\theta^+$ of Example~\ref{eg:monoidalkgraphs}. In particular, $\mathbb{F}_\theta^+$ embeds in its fundamental group.
\end{enumerate}
\end{examples}

\begin{example} \label{ex:bnex2}
Fix $n \ge 2$, and a permutation $\sigma \in \operatorname{Bij} ( [n] )$, the group of all bijections of the set $[n]$. Define $\theta : [n] \times [n] \to [n] \times [n]$ by $\theta (i,j) = (\sigma(i) ,j )$. This fits into the situation of Example~\ref{ex:bnex3}, so $\mathbb{F}_\theta^+$ embeds in its fundamental group.
\end{example}

\begin{dfn}
 Let $X$ be a nonempty set.  A map $R : X^2 \to X^2$ is a (set-theoretic) \emph{Yang--Baxter solution} if
 \[
 (R \times \text{id}_X)(\text{id}_X \times R)(R \times \text{id}_X) = (\text{id}_X \times R)(R \times \text{id}_X) (\text{id}_X \times R)
 \]
 as maps on $X^3$. For every permutation $\sigma$ of $X$ there is a Yang--Baxter solution
  $R$ given by $R(e, f) = (\sigma(f), e)$; such solutions are called
  \emph{permutation-type Yang--Baxter solutions}.
\end{dfn}

For details of the interplay between the Yang--Baxter equation and $k$-graphs see \cite{Y4}.

\begin{lem} \label{lem:ybex}
Fix a finite set $X$ and a Yang--Baxter solution $R : X^2 \to X^2$ on $X$. Fix $k \ge 2$. Let $\Lambda_{k, R}^0 = \{v\}$. For $i \le k$, let $\Lambda_{k, R}^{\varepsilon_i}=\{i\} \times X$. For $(i, e) \in \Lambda_{k, R}^{\varepsilon_i}$ and $(j, f) \in \Lambda_{k, R}^{\varepsilon_j}$ with $i < j$, set
\[
(i,e)(j,f) = (j, f')(i, e') \text{  if  } R(e,f) = ( f',e') .
\]
There is a unique $k$-graph $\Lambda_{k, R}$ with these edges and factorisation rules. If $R$ is a permutation type Yang--Baxter solution, then $i : \Lambda_{k, R} \to \Pi ( \Lambda_{k, R} )$ is injective.
\end{lem}
\begin{proof}
The first statement follows from \cite[\S 4.1]{Y4}. For the second statement we proceed by induction. For $k=2$ this follows from \cite[\S 4.1]{Y4}. Now suppose inductively that $\Lambda_{k-1, R}$ embeds in its fundamental groupoid. There is an automorphism $\alpha$ of $\Lambda_{k-1, R}$ such that $\alpha(i, e) = (i, \sigma(e))$ for all $i \le k-1$ and $e \in X$. For $e \in B^1_{|X|}$, let $\alpha_{e} := \alpha \in \Aut(\Lambda_{k-1, R})$. Corollary~\ref{thm:genaction} yields a $k$-graph $\Gamma = B_{|X|} \times_\alpha \Lambda_{k-1, R}$. Choose a bijection $\phi : B_{|X|}^1 \to \Lambda_{k, R}^{\varepsilon_1}$. Then there is an isomorphism $\Gamma \to \Lambda_{k, R}$ that agrees with $\phi$ on $B^1_{|X|} \subseteq \Gamma$ and takes each $(i, e) \in \Lambda^{\varepsilon_i}_{k-1, R} \subseteq \Gamma$ to $(i+1, e) \in \Lambda^{\varepsilon_{i+1}}_{k, R}$. Corollary~\ref{thm:genaction} implies that $\Gamma$ embeds in its fundamental groupoid, so $\Lambda_{k, R}$ does too.
\end{proof}

\begin{rmk}
For a long time the literature on $k$-graphs lacked concrete examples with $k \ge 3$ not obtained from lower-rank graphs via the constructions of Corollary~\ref{cor:specialcases}. Yang's important insight \cite{Y4} remedied this situation: every Yang--Baxter solution yields $k$-graphs for arbitrary $k$, typically not of the forms from Corollary~\ref{cor:specialcases}. In particular, Lemma~\ref{lem:ybex}, uses Yang's construction to see that every finite permutation $\sigma$ yields a $k$-graph that embeds in its fundamental groupoid for each $k \ge 1$. Taking $\sigma = \id$ yields cartesian-product $k$-graphs, but most other choices of $\sigma$ yield $k$-graphs that do not arise from the constructions of Corollary~\ref{cor:specialcases}.
\end{rmk}

\begin{rmk}\label{rmk:LV}
Work of Lawson and Vdovina also yields many embeddable $k$-graphs. A monoidal $k$-graph is \emph{rigid}
\cite[page~37]{LV} if whenever $e$ and $f$ are edges of different degrees, there are unique
edges $e', e'', f', f''$ such that $e'f = f'e$ and $ef'' = fe''$. Theorem~\ref{prop:prec} and
\cite[Theorem~11.14]{LV} combined imply that every rigid monoidal $k$-graph $\Lambda$ embeds in
$\Pi(\Lambda)$.
\end{rmk}

We finish the section by showing that a strongly connected $k$-graph $\Lambda$ embeds in $\Pi(\Lambda)$ whenever the submonoid of endomorphisms at any vertex embeds in a group.

\begin{prop}\label{prop:any vertex}
Let $\Lambda$ be a strongly connected $k$-graph and $H$ a group. Fix $v \in \Lambda^0$. If there exists an injective monoid homomorphism $c : v\Lambda v \to H$, then $i : \Lambda \to \Pi(\Lambda)$ is injective.
\end{prop}

\begin{proof}
The universal property of $\Pi(\Lambda)$ given in Definition~\ref{dfn:fundgpd} implies that there is a homomorphism $\tilde{c} : i(v)\Pi(\Lambda)i(v) \to H$ such that $\tilde{c} \circ i = c$.  Since $\Lambda$ is strongly connected, and since $\Pi(\Lambda)$ is a discrete groupoid, $\Pi(\Lambda)$ is isomorphic to $T(\Lambda^0) \times i(v)\Pi(\Lambda)i(v)$. Post-composing this isomorphism with $\id_{T(\Lambda^0)} \times \tilde{c}$ yields a groupoid homomorphism $q : \Pi(\Lambda) \to T(\Lambda^0) \times H$. Suppose that $q(i(\mu)) = q(i(\nu))$. Fix $\lambda \in v\Lambda r(\mu)$ and $\tau \in s(\mu)\Lambda v$. We have
\begin{align*}
	((v,v), c(\lambda\mu\tau))
		&= ((v,v), \tilde{c}(i(\lambda\mu\tau))
		= q(i(\lambda\mu\tau))\\
		&= q(i(\lambda\nu\tau))
		= ((v,v), \tilde{c}(i(\lambda\nu\tau))
		((v,v), c(\lambda\nu\tau)).
\end{align*}
Since $c$ is injective, $\lambda\mu\tau = \lambda\nu\tau$ and so $\mu = \nu$. Thus $q \circ i$, and therefore $i$, is injective.
\end{proof}

\begin{example}
Consider the $2$-graph $\Lambda$ below with relations $a_0e=fa_1,  \ a_1e=fa_0 ,  \ bf=eb$.
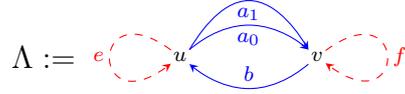
\[
\begin{tikzpicture}[xscale=0.9, yscale=0.65, >=stealth]
\node[inner sep=0.8pt] (u) at (0,0) {$\scriptstyle u$} edge[red, loop,->,in=225,out=135,looseness=15,dashed] node[auto,red,swap,inner sep=2pt] {$\scriptstyle e$} (u);
\node[inner sep=0.8pt] (v) at (2,0) {$\scriptstyle v$} edge[red, loop,->,in=-45,out=45,looseness=15,dashed] node[auto,red,right,inner sep=2pt] {$\scriptstyle f$} (v);
\draw[blue, ->] (u.north east)
    parabola[parabola height=0.5cm] (v.north west);
\node[blue, inner sep=0pt, below] at (1,0.5) {$\scriptstyle a_0$};
\draw[blue, ->] (v.south west)
    parabola[parabola height=-0.5cm] (u.south east);
\node[blue, inner sep=0pt, below] at (1,1) {$\scriptstyle a_1$};
\draw[blue, ->] (u.north east)
    parabola[parabola height=1cm] (v.north west);
\node[blue, inner sep=0pt, above] at (1,-0.5) {$\scriptstyle b$};
\node at (-2,0) {$\Lambda:=$};
\end{tikzpicture}
\]
Then $\Lambda$ is strongly connected. None of our results before Proposition~\ref{prop:any vertex} applies to show that $\Lambda$ embeds in $\Pi(\Lambda)$. Since $e b a_i = b f a_i = ba_{1-i} e$ for each $i$,  the monoid $u\Lambda u \subseteq \Lambda$ has presentation
\[
u \Lambda u = \langle e, b a_0 ,  b a_1 :  eba_i = ba_{1-i} e, i=0,1 \rangle,
\]
so is isomorphic to the semidirect product $\FF^+_2 \times_\alpha \NN$ for the action $\alpha$ that interchanges $\{ b a_0 ,  b a_1 \}$,
the generators of $\FF^+_2$. The action $\alpha$ extends uniquely to an action $\tilde{\alpha}$ of $\ZZ$ on $\FF_2$, and $u \Lambda u \cong \FF^+_2 \times_\alpha \NN$ embeds in $\FF_2 \times_{\tilde{\alpha}} \ZZ$. So Proposition~\ref{prop:any vertex} implies that $i : \Lambda \to \Pi(\Lambda)$ is injective.
\end{example}

\section{\texorpdfstring{$C^*$}{C*}-algebraic results} \label{sec:cstar}

Here we generalise \cite[Corollary~4.14]{KP1}, which says that the $C^*$-algebra of a connected row-finite $1$-graph is Rieffel--Morita equivalent to a crossed product of a commutative $C^*$-algebra by the fundamental group of the graph. The situation is much more complicated in higher dimensions.

Let $\Lambda$ be a connected row-finite source-free $k$-graph. Fix $v\in \Lambda^0$. By Theorem~\ref{thm:univisspg} (see \cite[Corollary~6.5]{PQR2}) there is a cocycle $\eta : \Lambda \to \pi_1 (\Lambda , v)$ such that the skew-product $ \pi_1 ( \Lambda , v ) \times_\eta \Lambda$ is isomorphic to the universal cover $\Sigma$ of $\Lambda$. It then follows from \cite[Theorem~5.7]{KP2} that $C^*(\Lambda)$ is Rieffel--Morita equivalent to $C^*(\Sigma)\rtimes \pi_1(\Lambda,v)$. Our main theorem describes the coefficient algebra $C^*(\Sigma)$ of this crossed product.

\begin{thm}\label{thm:mainC*result}
Let $\Lambda$ be a connected row-finite source-free $k$-graph, and let $\Sigma = \pi_1(\Lambda, v) \times_\eta \Lambda$ be as above so that $C^*(\Lambda)$ is Rieffel--Morita equivalent to $C^*(\Sigma) \rtimes \pi_1(\Lambda, v)$.
	\begin{enumerate}[(i)]
		\item The $C^*$-algebra $C^*(\Sigma)$ is AF.
		\item If $\Lambda$ embeds in its fundamental groupoid $\Pi(\Lambda)$, then $C^*(\Sigma)$ is type I$_0$ and its spectrum has a cover by zero-dimensional compact open Hausdorff subsets.
		\item If $\Lambda$ embeds in its fundamental groupoid $\Pi(\Lambda)$ and $\Sigma^{\mathbb{N}\mathbf{1}}$ is simply connected, then $C^*(\Sigma)$ is Rieffel--Morita equivalent to a commutative $C^*$-algebra.
	\end{enumerate}
\end{thm}

We use the next two results to prove parts (i)~and~(ii) of Theorem~\ref{thm:mainC*result}.

\begin{prop}\label{thm:univisAF}
	Let $\Gamma$ be a row-finite source-free $k$-graph. If $\Gamma$ is simply connected, then there is a map $f: \Gamma^0 \to  \ZZ^k$  such that $d(\lambda) = f(s(\lambda)) -  f(r(\lambda))$
	for all $\lambda \in \Gamma$. Moreover, $C^*(\Gamma)$ is AF.
\end{prop}

\begin{proof}
	Since $d : \Gamma \to \ZZ^k$ is a cocycle, Lemma~\ref{lem:characterise simply connected} ensures the existence of $f$. Now \cite[Lemma~5.4]{KP2} implies that  $C^*(\Gamma)$ is AF.
\end{proof}

\begin{example}
Let $\Lambda$ be the $2$-graph of Example~\ref{ex:pqrex1} (see \cite[Example~7.1]{PQR1}).
Recall that $\Lambda$ does not embed in its fundamental groupoid and the universal cover $\Sigma = \ZZ^2 \times_d \Lambda$ is simply connected.
We claim that $C^*(\Sigma)$ is Rieffel--Morita equivalent to the UHF algebra $M_{6^{\infty}}$ (in fact, $C^*(\Sigma) \cong M_{6^{\infty}} \otimes \mathcal{K}$).
For each $n \in \mathbb{N}$, set $v_n := n\bf{1}$ and observe that as in the proof of \cite[Lemma~5.4]{KP2},  
$A_n := C^*(\{ s_\lambda : s(\lambda) = v_n \})  \cong \mathcal{K}(\ell^2(s^{-1}(v_n)))$.
Moreover, for all $n$, $A_n \subset A_{n+1}$ and the multiplicity of the embedding is 6 (since $|v_n\Lambda v_{n+1}| = 6$).
Since the sequence of $v_n$'s is cofinal in $\ZZ^2$ we have
\[
C^*(\Sigma) \cong \varinjlim A_n \cong  \varinjlim \mathcal{K}(\ell^2(s^{-1}(v_n))).
\]
Hence,  $C^*(\Sigma)$ is Rieffel--Morita equivalent to the UHF algebra $M_{6^{\infty}}$ as claimed.
\end{example}

\begin{prop}\label{prop:fell}
	Let $\Gamma$ be a row-finite source-free $k$-graph.	If $\Gamma$ is singly connected, then, for each $v \in \Gamma^0$, the corner $s_vC^*(\Gamma)s_v$ is an abelian $C^*$-algebra isomorphic to $C(Z(v))$. Moreover, $C^*(\Gamma)$ is type I$_0$, and $\operatorname{Prim} C^*(\Gamma)$ admits a cover by zero-dimensional compact open Hausdorff sets.
\end{prop}
\begin{proof}
	Fix $x, y \in \Gamma^\infty$ such that $x(0) = y(0)$, and $p,q\in\NN^k$. We claim that if $\sigma^p(x) = \sigma^q(y)$, then  $p =q$ and $x = y$. To see this, suppose that $\sigma^p(x) = \sigma^q(y)$.  Then $\sigma^p(x)(0) = \sigma^q(y)(0)$. Let $u := x(0) = y(0)$ and $v := \sigma^p(x)(0) = \sigma^q(y)(0)$. Then $x(0, p), y(0, q) \in u{\Gamma}v$.
	Since $\Gamma$ is singly connected, $x(0, p) = y(0, q)$.  Hence,
	$x = x(0, p)\sigma^p(x) = y(0, q)\sigma^q(y)(0) = y$, and the claim holds.
	
	Now recall from \cite{KP2} that $\Gamma^\infty = \Gg_\Gamma^0$, and that for $v \in \Gamma^0$, we have
	\[
	\Gg_\Gamma|_{Z(v)} := \{ \gamma \in \Gg_\Gamma : s(\gamma), r(\gamma) \in Z(v) \},
	\]
     $s_v = \chi_{Z(v)} \subset C_0(\Gg_\Gamma^0)$, and $s_vC^*(\Gamma)s_v \cong C^*( \Gg_\Gamma|_{Z(v)})$.	By the first paragraph, $\Gg_\Gamma|_{Z(v)}^0 \cong Z(v)$. Hence $C^*( \Gg_\Gamma|_{Z(v)}) \cong C(Z(v))$.	So for each $v \in \Gamma^0$, the ideal $I_v$ generated by $s_v$ is Rieffel--Morita equivalent to the abelian $C^*$-algebra $C(Z(v))$. Since $C^*(\Gamma)$ is generated by the ideals $I_v$, $C^*(\Gamma)$ is type I$_0$.

    By definition of the hull-kernel topology, the ideals $I_v$ yield a cover of $\operatorname{Prim}(C^*(\Gamma))$ by open sets $\widehat{I}_v \cong \operatorname{Prim}(I_v)$. Since each $I_v$ is Rieffel--Morita equivalent to $C(Z(v))$, each $\operatorname{Prim}(I_v) \cong Z(v)$ is a zero-dimensional compact open Hausdorff subspace of $\operatorname{Prim}(C^*(\Gamma))$.
\end{proof}

\begin{proof}[Proof of Theorem~\ref{thm:mainC*result}~(i)~and~(ii)]
Proposition~\ref{thm:univisAF} for $\Gamma=\Sigma$ gives~(i). If $\Lambda\to\Pi(\Lambda)$ is injective then so is $\Sigma\to\Pi(\Sigma)$ by Proposition~\ref{prop:lift}. Since $\Sigma$ is simply connected, Proposition~\ref{prop:1cnc} implies that $\Sigma$ is singly connected; so Proposition~\ref{prop:fell} for $\Gamma=\Sigma$ gives~(ii).
\end{proof}

To prove Theorem~\ref{thm:mainC*result}(iii), we will argue that the spectrum $C^*(\Sigma)^{\wedge}$ is Hausdorff: then Theorem~\ref{thm:mainC*result}(2) shows that $C^*(\Sigma)$ is Rieffel--Morita a continuous-trace $C^*$-algebra with totally disconnected spectrum, and the Dixmier--Douady theorem will show that $C^*(\Sigma)$ is Rieffel--Morita equivalent to $C_0(C^*(\Sigma)^{\wedge})$. We will argue that $C_0(C^*(\Sigma)^{\wedge}) \cong \Sigma^{\NN\mathbf{1}}/\Gg_{\Sigma^{\NN\mathbf{1}}}$, and use the additional hypothesis that $\Sigma^{\NN\mathbf{1}}$ is simply connected to prove Theorem~\ref{thm:mainC*result}(3). We do not know whether this additional hypothesis is automatic; certainly, even for $1$-graphs, being singly connected does not guarantee that the associated $C^*$-algebra has Hausdorff spectrum:

\begin{example}\label{eg:nonHausdorff}
	Let $E$ be the directed graph (pictured below) such that
	\begin{itemize}
		\item $E^0 = \{u_n, v_n : n \in \ZZ\} \cup \{w_{n, i} : n \in \ZZ \text{ and }i \ge 0\}$, and
		\item $E^1 = \{e_n, f_n, g_n, h_n : n \in \ZZ\} \cup \{k_{n, i} : n \in \ZZ\text{ and }i \ge 0\}$,
	\end{itemize}
    and such that for $n \in \ZZ$ and $i \ge 0$,
	\begin{align*}
		r(e_n) &= s(e_{n-1}) = r(g_n) = u_n, &  r(f_n) &= s(f_{n-1}) = r(h_n) = v_n,\\
		s(g_n) &= s(h_n) = r(k_{n,0}) = w_{n,0}, & s(k_{n,i}) &= r(k_{n, i+1}) = w_{n, i+1}.
	\end{align*}
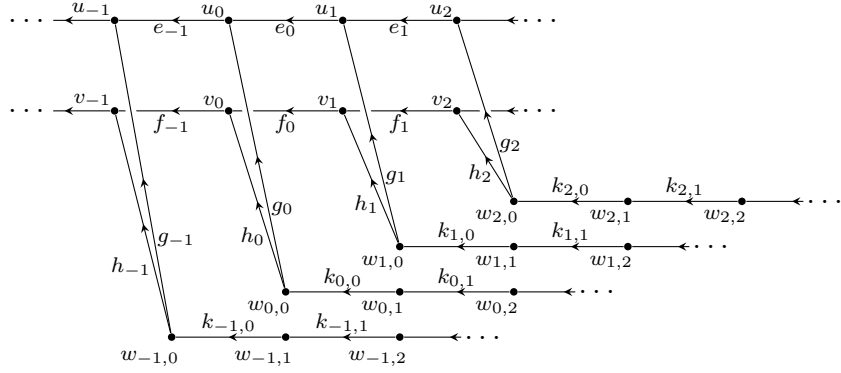
\[
\begin{tikzpicture}[yscale=1.2, xscale=1.5, z={(0.5,-1)}, >=stealth]
    \node[inner sep=1pt] (uldots) at (-1.75,1,0) {\dots};
    \node[inner sep=1pt] (vldots) at (-1.75,0,0) {\dots};
    \foreach \x in {-1,0,1,2} {
        \node[fill=black, circle, inner sep=1pt] (u\x) at (\x,1,0) {};
        \node[anchor=south east, inner sep=1pt] at (u\x) {\tiny$u_{\x}$};
        \node[fill=black, circle, inner sep=1pt] (v\x) at (\x,0,0) {};
        \node[anchor=south east, inner sep=1pt] at (v\x) {\tiny$v_{\x}$};
    }
    \node[inner sep=1pt] (urdots) at (2.75,1,0) {\dots};
    \node[inner sep=1pt] (vrdots) at (2.75,0,0) {\dots};
    \foreach \x in {0,1,2} {
            \node[fill=black, circle, inner sep=1pt] (w-1\x) at (\x-1,-1.5,1) {};
            \node[anchor=north east, circle, inner sep=0pt] at (w-1\x) {\tiny$w_{-1,\x}$};
    }
    \node[inner sep=1pt] (w-1dots) at (1.75,-1.5,1) {\dots};
    \foreach \x in {0,1,2} {
       \node[fill=black, circle, inner sep=1pt] (w0\x) at (\x,-1,1) {};
       \node[anchor=north east, circle, inner sep=0pt] at (w0\x) {\tiny$w_{0,\x}$};
    }
    \node[inner sep=1pt] (w0dots) at (2.75,-1,1) {\dots};
    \foreach \x in {0,1,2} {
            \node[fill=black, circle, inner sep=1pt] (w1\x) at (\x+1,-0.5,1) {};
            \node[anchor=north east, circle, inner sep=0pt] at (w1\x) {\tiny$w_{1,\x}$};
    }
    \node[inner sep=1pt] (w1dots) at (3.75,-0.5,1) {\dots};
    \foreach \x in {0,1,2} {
            \node[fill=black, circle, inner sep=1pt] (w2\x) at (\x+2,0,1) {};
            \node[anchor=north east, circle, inner sep=0pt] at (w2\x) {\tiny$w_{2,\x}$};
    }
    \node[inner sep=1pt] (w2dots) at (4.75,0,1) {\dots};
    \begin{scope}[decoration={markings, mark=at position 0.2 with {\arrow{>}}}]
        \draw[postaction=decorate] (urdots) to (u2);
        \draw[postaction=decorate] (vrdots) to (v2);
        \draw[postaction=decorate] (w-1dots) to (w-12);
        \draw[postaction=decorate] (w0dots) to (w02);
        \draw[postaction=decorate] (w1dots) to (w12);
        \draw[postaction=decorate] (w2dots) to (w22);
    \end{scope}
    \begin{scope}[decoration={markings, mark=at position 0.8 with {\arrow{>}}}]
        \draw[postaction=decorate] (u-1) to (uldots);
        \draw[postaction=decorate] (v-1) to (vldots);
    \end{scope}
    \begin{scope}[decoration={markings, mark=at position 0.5 with {\arrow{>}}}]
    \foreach \x/\w in {-1/0,0/1,1/2} {
        \draw[postaction={decorate}] (u\w) to node[below, pos=0.5, inner sep=1pt]{\tiny$e_{\x}$} (u\x);
        \draw[postaction={decorate}] (v\w) to node[below, pos=0.5, inner sep=1pt]{\tiny$f_{\x}$} (v\x);
    }
    \draw[line width=4pt, white] (w-10)--(u-1);
    \draw[line width=4pt, white] (w-10)--(v-1);
    \draw[postaction=decorate] (w-10) to node[right, pos=0.3, inner sep=1pt] {\tiny$g_{-1}$} (u-1);
    \draw[postaction=decorate] (w-10) to node[left, pos=0.3, inner sep=1pt] {\tiny$h_{-1}$} (v-1);
    \draw[postaction=decorate] (w-11) to node[above, pos=0.5, inner sep=1pt] {\tiny$k_{-1,0}$} (w-10);
    \draw[postaction=decorate] (w-12) to node[above, pos=0.5, inner sep=1pt] {\tiny$k_{-1,1}$} (w-11);
    \draw[line width=4pt, white] (w00)--(u0);
    \draw[line width=4pt, white] (w00)--(v0);
    \draw[postaction=decorate] (w00) to node[right, pos=0.3, inner sep=1pt] {\tiny$g_0$} (u0);
    \draw[postaction=decorate] (w00) to node[left, pos=0.3, inner sep=1pt] {\tiny$h_0$} (v0);
    \draw[postaction=decorate] (w01) to node[above, pos=0.5, inner sep=1pt] {\tiny$k_{0,0}$} (w00);
    \draw[postaction=decorate] (w02) to node[above, pos=0.5, inner sep=1pt] {\tiny$k_{0,1}$} (w01);
    \draw[line width=4pt, white] (w10)--(u1);
    \draw[line width=4pt, white] (w10)--(v1);
    \draw[postaction=decorate] (w10) to node[right, pos=0.3, inner sep=1pt] {\tiny$g_1$} (u1);
    \draw[postaction=decorate] (w10) to node[left, pos=0.3, inner sep=1pt] {\tiny$h_1$} (v1);
    \draw[postaction=decorate] (w11) to node[above, pos=0.5, inner sep=1pt] {\tiny$k_{1,0}$} (w10);
    \draw[postaction=decorate] (w12) to node[above, pos=0.5, inner sep=1pt] {\tiny$k_{1,1}$} (w11);
    \draw[line width=4pt, white] (w20)--(u2);
    \draw[line width=4pt, white] (w20)--(v2);
    \draw[postaction=decorate] (w20) to node[right, pos=0.3, inner sep=1pt] {\tiny$g_2$}(u2);
    \draw[postaction=decorate] (w20) to node[left, pos=0.3, inner sep=1pt] {\tiny$h_2$} (v2);
    \draw[postaction=decorate] (w21) to node[above, pos=0.5, inner sep=1pt] {\tiny$k_{2,0}$} (w20);
    \draw[postaction=decorate] (w22) to node[above, pos=0.5, inner sep=1pt] {\tiny$k_{2,1}$} (w21);
    \end{scope}
\end{tikzpicture}
\]

This graph $E$ is singly connected. Define $x,y \in E^\infty$ by $x = e_0 e_1 e_2 \cdots$ and $y = f_0 f_1 f_2\cdots$. Then $[x] \not= [y]$  in $\Gg_E^{(0)}/\Gg_E$ . We claim that they cannot be separated by disjoint open sets.  To see this, for $n \in \ZZ$, let $z_n = k_{n,0} k_{n, 1} k_{n, 2} \cdots$. We will show that $[z_n] \to [x]$ and $z_n \to [y]$ as $n \to \infty$. By symmetry, we just have to show that $[z_n] \to [x]$. For this, just note that $[z_n] = [e_0 e_1 \cdots e_{n-1} g_n z_n]$, and we have $\lim_{n \to \infty} e_0 e_1 \dots e_{n-1} g_n z_n = e_0 e_1 e_2 \cdots = x$.

We have $C^*(E) \cong C^*(\Gg_E)$ by \cite[Proposition 4.1]{KPRR}. Since $C^*(E)$ is type I$_0$, its spectrum is homeomorphic, by \cite[Corollary~4.2]{OrloffClark}, to the orbit space $\Gg_E^{(0)}/\Gg_E$ of $\Gg_E$, which we just saw is not Hausdorff. Note that $E$ is not simply connected (for example $e_0g_1h_1^{-1}f_0^{-1}h_0g_0^{-1} \in \Pi(E^*)^{u_0}_{u_0} \setminus \{u_0\}$).
\end{example}

Example~\ref{eg:nonHausdorff} suggests a Hausdorffness criterion (Lemma~\ref{lem:paths version Hausdorff}). As this criterion is not easy to check, in Theorem~\ref{thm:vertices version Hausdorff}, we specialise to singly connected $k$-graphs and recast it in terms of the following relation on vertices, which permeates analyses of ideals of $k$-graph $C^*$-algebras \cite{RSY1}.

\begin{nota}
For $\Gamma$ a $k$-graph, we define a relation $\le$ on $\Gamma^0$ by $v \le w$ if and only if $v \Gamma w \not= \emptyset$.
\end{nota}

\begin{lem}\label{lem:paths version Hausdorff}
Let $\Gamma$ be a row-finite source-free $k$-graph, and let $\Gg_\Gamma$ be its $k$-graph groupoid. The orbit space $\Gamma^\infty/\Gg_\Gamma$ is Hausdorff if and only if for every pair of infinite paths $x,y \in \Gamma^\infty$ such that $[x]\not=[y]$, there exists $N \in \NN^k$ such that the vertices $x(N)$ on $x$ and $y(N)$ on $y$ have no common upper bound with respect to $\le$, in the sense that $s(\mu) \not= s(\nu)$ for all $\mu \in x(N)\Gamma$ and $\nu \in y(N)\Gamma$.
\end{lem}
\begin{proof}
We have $[x] = [y]$ if and only if $\sigma^m(x) = \sigma^n(y)$ for some $m,n$. So it suffices to fix $x,y$ such that $\sigma^m(x) \not= \sigma^n(y)$ for all $m,n$, and show that $[x]$ and $[y]$ can be separated if and only if there exists $N$ as in the statement. Suppose that there is no such $N$. For each $N \in \NN^k$, choose $\mu_N \in x(N)\Gamma$ and $\nu_N \in y(N)\Gamma$ with $s(\mu_N) = s(\nu_N)$, and $z_N \in s(\mu_N)\Gamma^\infty$. Then $x(0,N)\mu_Nz_N \to x$ and $y(0, N)\nu_N z_N \to y$. Since each $[x(0,N)\mu_Nz_N] = [z_N] = [y(0,N)\nu_Nz_N]$, this forces $[z_N] \to [x]$ and $[z_N] \to [y]$. Now suppose that there exists $N$ as in the statement. Then $q(Z(x(0, N)))$ and $q(Z(y(0,N)))$ are disjoint open neighbourhoods of $[x]$ and $[y]$ in $\Gamma^\infty/\Gg_\Gamma$.
\end{proof}

Recall that a $\emph{filter}$ for a partially ordered set $(X, \preceq)$ is a nonempty subset $F \subseteq X$ such that
\begin{itemize}
\item[(a)] for all $u,v \in F$ there exists $w \in F$ such that $u,v \preceq w$;
\item[(b)] if $v \in F$ and $u \preceq v$, then $u \in F$.
\end{itemize}
A filter $F$ for $\preceq$ is called an \emph{ultrafilter} if
\begin{itemize}
\item[(c)] $F$ is not properly contained in any other filter $F'$ for $(X, \preceq)$.
\end{itemize}

If $\Gamma$ is singly connected, then $\le$ is a partial order on $\Gamma^0$. We show that elements of $\Gamma^\infty/\Gg_\Gamma$ correspond with ultrafilters for $(\Gamma^0, \le)$, and use this to characterise Hausdorffness of $\Gamma^\infty/\Gg_\Gamma$.

\begin{thm}\label{thm:vertices version Hausdorff}
Let $\Gamma$ be a singly connected row-finite source-free $k$-graph. Then the ultrafilters for $(\Gamma^0, \le)$ are exactly the sets $[x]^0 := \{r(y) : y \in [x]\}$ indexed elements $x \in \Gamma^\infty$. Moreover,  $\Gamma^\infty/\Gg_\Gamma$ is Hausdorff if and only if for every pair $U, V$ of distinct ultrafilters of $(\Gamma^0, \le)$ there is a pair $u \in U$ and $v \in V$ with no common upper bound with respect to $\le$.
\end{thm}
\begin{proof}
For the first statement, first fix $x \in \Gamma^\infty$. If $v_1, v_2 \in [x]^0$, then $v_1 = r(\alpha \sigma^m(x))$ and $v_2 = r(\beta\sigma^n(x))$ for some $\alpha,\beta, m,n$ and then $w = r(\sigma^{m+n}(x)) \in [x]^0$ satisfies $v_i \Gamma w \not= \emptyset$ by definition; so $[x]^0$ satisfies~(a). If $w \in [x]^0$ and $v \in \Gamma^0$ satisfy $v \le w$, say $\alpha \in v \Gamma w$, then since $w \in [x]^0$ we have $w = r(\beta \sigma^{n}(x))$ for some $\beta, n$ and so $v = r(\alpha\beta \sigma^n(x)) \in [x]^0$; so $[x]^0$ satisfies~(b). Suppose that $F$ is a filter for $(\Gamma^0, \le)$ containing $[x]^0$. Fix $v \in F$; we must show that $v \in [x]^0$. Since $v, x(0) \in F$ there exists $w \in F$ with $v \le w$ and $x(0) \le w$, and by~(b), if $w \in [x]^0$ then $v \in [x]^0$; so we just have to show that $w \in [x]^0$. Fix $\alpha \in x(0)\Gamma w$. Then $r(\sigma^{d(\alpha)}(x)) \in [x]^0 \subseteq F$. So there exists $w' \in F$ such that $r(\sigma^{d(\alpha)}(x)), w \le w'$; say $\mu \in r(\sigma^{d(\alpha)}(x)) \Gamma w'$ and $\nu \in w \Gamma w'$. So $\alpha \nu$ and $x(0, d(\alpha))\mu$ both belong to $x(0) \Gamma w'$. Since $\Gamma$ is singly connected, this forces $\alpha\nu = x(0, d(\alpha))\mu$, so the factorisation property forces $x(0, d(\alpha)) = \alpha$; hence $w = s(\alpha) = s(x(0,d(\alpha))) \in [x]^0$.

Now fix an ultrafilter $F$ for $(\Gamma^0, \le)$. Enumerate $F = (v_1, v_2, \dots)$, put $w_1 = v_1$ and inductively use~(a) to choose $w_{i+1} \in F$ such that $v_{i+1}, w_i \le w_{i+1}$. So $(w_i)_i$ is an increasing sequence such that every $v \in F$ satisfies $v \le w_i$ for some $i$. For each $i$, use that $w_i \le w_{i+1}$ to fix $\alpha_i \in w_i \Gamma w_{i+1}$, let $\mu_i := \alpha_1\dots\alpha_i$, and choose $y_i \in Z(\mu_i) \subseteq \Gamma^\infty$. Since $(y_i)_i$ belongs to the compact set $Z(v_1)$ it has a convergent subsequence $y_{i_l} \to y \in Z(v_1)$. We claim that $F = [y]^0$. By~(c) it suffices to show that $F \subseteq [y]^0$. So fix $v \in F$. Then $v = v_m \le w_m$ for some $m \in \NN$. Choose $l$ so that $i_l \ge m$. For $l' \ge l$ we have $y_{i_{l'}} \in Z(\mu_{i_{l'}}) \subseteq Z(\mu_{i_l})$. Hence $y \in Z(\mu_{i_l})$. So $w_{i_l} = s(\mu_{i_l}) = r(\sigma^{d(\mu_{i_l})}(y)) \in [y]^0$. By choice $(w_j)_j$, and $l$, we have $v \le w_m \le w_{i_l}$. So~(b) gives $v \in [y]^0$. This proves the first statement.

For the second statement, by Lemma~\ref{lem:paths version Hausdorff}, it suffices to show that for all $x,y \in \Gamma^\infty$ with $[x] \not= [y]$, there exists $N$ such that $s(x(N)\Gamma) \cap s(y(N)\Gamma) = \emptyset$ if and only if, for all pairs $U \not= V$ of ultrafilters of $(\Gamma^0, \le)$ there exist $u \in U$ and $v \in V$ with no common upper bound with respect to $\le$.

First suppose that for every pair $x,y \in \Gamma^\infty$ there exists $N$ such that $s(\mu) \not= s(\nu)$ for every $\mu \in x(N)\Gamma$ and $\nu \in y(N)\Gamma$. Fix ultrafilters $U \not= V$, and fix $x,y \in \Gamma^\infty$ with $U = [x]^0$ and $V = [y]^0$. Fix $N$ such that $s(\mu) \not= s(\nu)$ for every $\mu \in x(N)\Gamma$ and $\nu \in y(N)\Gamma$. Then $u = x(N) \in U$ and $v = y(N) \in V$ have no common upper bound. Now suppose that for every pair $U \not= V$ of ultrafilters, there exist $u \in U$ and $v \in V$ with no common upper bound. Fix $x,y \in \Gamma^\infty$ with $U = [x]^0$ and $V = [y]^0$. Fix $u \in U$ and $v \in V$ with no common upper bound. Fix $x' \in [x]$ and $y' \in [y]$ with $r(x') = u$ and $r(y') = v$, and $m, m'$ and $n, n'$ such that $\sigma^m(x) = \sigma^{m'}(x')$ and $\sigma^n(y) = \sigma^{n'}(y')$. Fix $N \ge m,n$. Then $u \Gamma x(N) \not= \emptyset$ and $v \Gamma y(N) \not= \emptyset$. Since $u,v$ have no common upper bound, nor do $x(N)$ and $y(N)$; so $s(\mu) \not= s(\nu)$ for all $\mu \in x(N)\Gamma$ and $\nu \in y(N)\Gamma$.
\end{proof}

\begin{rmk}\label{rmk:tree ends Hausdorff}
Lemma~\ref{lem:paths version Hausdorff} gels with \cite[Proposition~4.3]{KP1}: if $E$ is a simply connected row-finite source-free directed graph, then $E^\infty/\Gg_E$ is Hausdorff. We prove the contrapositive. Suppose that $E^\infty/\Gg_E$ is not Hausdorff. Since $E$ is a $1$-graph, $i : E^* \to \Pi(E^*)$ is injective. Corollary~\ref{lem:paths version Hausdorff} gives $x,y \in E^\infty$ such that $\sigma^m(x) \not= \sigma^n(y)$ for all $m,n$, and $\mu_N \in x(N)E^*$ and $\nu_N \in y(N) E^*$ such that $s(\mu_N) = s(\nu_N) =: w_N$. We first claim that there exists $N_0$ such that $x(n) \not= y(m)$ for all $m,n \ge N_0$. To see this, suppose that there are increasing sequences $(n_i), (m_i)$ such that $x(n_i) = y(m_i)$ for all $i$. Since $E$ is singly connected, $x(n_i, n_{i+1}) = y(m_i, m_{i+1})$ for all $i$; so $\sigma^{n_0}(x) = \sigma^{n_0}(y)$, a contradiction. So by replacing $x, y$ with $\sigma^{N_0}(x)$ and $\sigma^{N_0}(y)$, we may assume that $x(m) \not= y(n)$ for all $m,n$.

Hence each $w_N$ is exactly one of $x, y$; without loss of generality, $w_0$ is not on $x$. Let $\alpha = x(0, |\mu_0|)$ and $\beta = y(0, |\nu_0|)$. Then $\mu_0^{-1}\alpha\mu_{|\mu_0|}\nu_{|\nu_0|}^{-1}\beta^{-1}\nu_0 \in (\Pi (E))^{w_0}_{w_0}$. We show that $\mu_0^{-1}\alpha\mu_{|\mu_0|}\nu_{|\nu_0|}^{-1}\beta^{-1}\nu_0 \not= w_0$. Since $r(\mu_{|\mu_0|}) = x(|\mu_0|) \not= y(|\nu_0|) = r(\nu_{|\nu_0|})$, in reduced form, $\mu_{|\mu_0|} \nu_{|\nu_0|}^{-1} = e \gamma\lambda^{-1}$ where $e \in E^1$ is the first edge of $\mu_{|\mu_0|}$. Similarly, since $w_0 = s(\nu_0) \not= x(|\mu_0|) = s(\alpha)$, in reduced form, $\mu_0^{-1}\alpha = \zeta^{-1}\eta f$ where $f$ is the last edge of $\alpha$ and $\zeta,\eta \in E^*$. So in reduced form $\mu_0^{-1}\alpha \mu_{|\mu_0|}\nu_{|\nu_0|}^{-1} = \zeta^{-1} \eta f e \gamma\lambda^{-1}$. In particular, the word $fe$ appears in the reduced form of $\mu_0^{-1}\alpha\mu_{|\mu_0|}\nu_{|\nu_0|}^{-1}\beta^{-1}\nu_0$, so this is a nontrivial element of $(\Pi (E))^{w_0}_{w_0}$. Hence $E$ is not simply connected.
\end{rmk}

\begin{rmk}
The argument of the preceding remark does not go through for $k$-graphs because there is no canonical reduced form for elements of the fundamental groupoid of a $k$-graph.
\end{rmk}

For $1$-graphs $E^*$ we can use Remark~\ref{rmk:tree ends Hausdorff} to check Hausdorffness of $E^\infty/\Gg_E \cong C^*(E^*)^{\wedge}$. So it helps to relate Hausdorffness of the orbit space of a $k$-graph to that of a natural sub-$1$-graph

\begin{prop}\label{prp:Hausdorffness}
Let $\Gamma$ be a row-finite source-free $k$-graph. Suppose that $\Gamma$ is simply connected. Then $\Gamma^\infty/\Gg_\Gamma$ is Hausdorff if and only if $(\Gamma^{\NN\mathbf{1}})^\infty/\Gg_{\Gamma^{\NN\mathbf{1}}}$ is Hausdorff.
\end{prop}

To prove this, we show that $\Gamma^\infty/\Gg_\Gamma$ is homeomorphic to a clopen subset of $(\Gamma^{\NN\mathbf{1}})^\infty/\Gg_{\Gamma^{\NN\mathbf{1}}}$.

\begin{lem}\label{lem:orbit spaces homeomorphic}
Let $\Gamma$ be a row-finite source-free $k$-graph. Suppose that $\Gamma$ is simply connected. Let $f : \Gamma^0 \to \ZZ^k$ be a function such that $d(\lambda) = f(s(\lambda)) - f(r(\lambda))$ for all $\lambda \in \Gamma$ as in Proposition~\ref{thm:univisAF}. Let $E$ be the directed graph such that $E^0 = f^{-1}(\ZZ\mathbf{1})$ and $E^1 = E^0 \Gamma^{\mathbf{1}}$. Let $j : E^\infty \to \Gamma^\infty$ be the map such that $j(x)$ is the unique infinite path such that $j(x)(0, n\cdot\mathbf{1}) = x_1 x_2\cdots x_n$ for all $n \in \NN$ (see \cite[Remark~2.2]{KP2}). Then $j$ descends to a homeomorphism $\widetilde{j} : E^\infty/\Gg_E \to \Gamma^\infty/\Gg_\Gamma$.
\end{lem}
\begin{proof}
The map $j$ restricts to a homeomorphism $vE^\infty \to v\Gamma^\infty$ for each $v \in E^0$, so is continuous.

We claim that if $x,y \in E^\infty$ then $j(x) \sim_{\Gg_\Gamma} j(y)$ if and only if $x \sim_{\Gg_E} y$. To see this, fix $x,y \in E^\infty$. Then $j(x) \sim_{\Gg_\Gamma} j(y)$ if and only if there exist $m,n \in \NN^k$ such that $\sigma^m(j(x)) = \sigma^n(j(y))$. Since $f(r(\sigma^m(j(x))) = f(r(x)) + m$ for all $m \in \NN^k$ and similarly for $y$, and since $f(r(x)), f(r(y)) \in \ZZ\mathbf{1}$, we deduce that $j(x) \sim_{\Gg_\Gamma} j(y)$ if and only if there exist $m,n \in \NN^k$ such that $\sigma^m(j(x)) = \sigma^n(j(y))$ and $m - n \in \ZZ\mathbf{1}$. Since $m - n \in \ZZ\mathbf{1}$ if and only if there exists $p \in \NN^k$ such that $m+p, n+p \in \NN\mathbf{1}$, we deduce that $j(x) \sim_{\Gg_\Gamma} j(y)$ if and only if $\sigma^{a\mathbf{1}}(j(x)) = \sigma^{b\mathbf{1}}(j(y))$ for some $a,b \in \NN$; that is, if and only if $\sigma^a(x) = \sigma^b(y)$ for some $a,b \in \NN$. Hence $j(x) \sim_{\Gg_\Gamma} j(y)$ if and only if $x \sim_{\Gg_E} y$. It follows that $j$ descends to a continuous function $\widetilde{j} : E^\infty/\Gg_E \to \Gamma^\infty/\Gg_\Gamma$.

Fix $p : \Gamma^0 \to \NN^k$ satisfying $f(v) + p(v) \in \ZZ\mathbf{1}$ for all $v$. For $x \in \Gamma^\infty$ and $j \in \NN$, let $\tilde{x}_j := \sigma^{p(r(x))}((j-1)\mathbf{1}, j\mathbf{1}) \in E^1$, and define $h(x) := \tilde{x}_1  \tilde{x}_2 \cdots  \tilde{x}_n \cdots \in E^\infty$. As $x \mapsto p(r(x))$ is locally constant, $h$ is continuous. Since $\sigma^{p(x)}(x) \sim_{\Gg_\Gamma} x$ for all $x$, the claim above shows that $x \sim_{\Gg_\Gamma} y$ if and only if $h(x) \sim_{\Gg_E} h(y)$, so $h$ descends to a continuous function $\tilde{h} : \Gamma^\infty/\Gg_\Gamma$ to $E^\infty/\Gg_E$. It is routine to check that $\tilde{h}$ and $\tilde{j}$ are mutually inverse: $h \circ j = \operatorname{id}_{E^\infty}$, and $[j \circ h(x)] = [\sigma^{p(x)}(x)] = [x]$ for all $x \in \Gamma^\infty$. In particular, $j$ descends to a homeomorphism as claimed.
\end{proof}

\begin{proof}[Proof of Proposition~\ref{prp:Hausdorffness}]
Resume the notation of Lemma~\ref{lem:orbit spaces homeomorphic}. It suffices for us to show that $(\Gamma^{\NN\mathbf{1}})^\infty/\Gg_{\Gamma^{\NN\mathbf{1}}}$ is Hausdorff if and only if $E^\infty/\Gg_E$ is Hausdorff.

For $p \in \ZZ^k$, let $V_p := f^{-1}(p + \ZZ\mathbf{1}) \subseteq \Gamma^0$ (so $V_0$ is $V$ in Lemma~\ref{lem:orbit spaces homeomorphic}). If $p-q \not \in \ZZ\mathbf{1}$, then $V_p\Gamma^{\NN\mathbf{1}} V_q = \emptyset$. So if $x \in V_p \Gamma^{\NN\mathbf{1}}$ and $y \in V_q \Gamma^{\NN\mathbf{1}}$, then $\sigma^{a\mathbf{1}}(x) \not= \sigma^{b\mathbf{1}}(y)$ for all $a,b \in \NN$ and hence $[x]_{\Gg_{\Gamma^{\NN\mathbf{1}}}} \not= [y]_{\Gg_{\Gamma^{\NN\mathbf{1}}}}$. Hence the sets $\big\{V_p (\Gamma^{\NN\mathbf{1}})^{\infty} : p \in \ZZ^{k-1} \times \{0\}\big\}$ have mutually disjoint open images in $(\Gamma^{\NN\mathbf{1}})^\infty/\Gg_{\Gamma^{\NN\mathbf{1}}}$. So it suffices to show that each of these images is Hausdorff.

Let $q : (\Gamma^{\NN\mathbf{1}})^\infty \to (\Gamma^{\NN\mathbf{1}})^\infty/\Gg_{\Gamma^{\NN\mathbf{1}}}$ be the quotient map. By assumption $E^\infty/\Gg_E = q(V_0 (\Gamma^{\NN\mathbf{1}})^{\infty})$ is Hausdorff, so it suffices to fix $p \in \ZZ^{k-1} \setminus\{0\}$ and show that $q(V_p\Gamma^{\NN\mathbf{1}})^\infty) \cong q(V_0\Gamma^{\NN\mathbf{1}})^\infty)$.

Since $V_p = V_{p + a\mathbf{1}}$ for all $a \in \NN$, we may assume that $p \ge 0$. Fix $n \in \NN^k$ such that $p + n \in \ZZ\mathbf{1}$. Then $\sigma^p : V_0 \Gamma^\infty \to V_p\Gamma^\infty$, and $\sigma^n : V_p \Gamma^\infty \to V_{p+n}\Gamma^\infty = V_0\Gamma^\infty$ are continuous. Using \cite[Remark~2.2]{KP2}, we can identify $V_p \Gamma^\infty$ with $V_p(\Gamma^{\NN\mathbf{1}})^\infty$ and $V_0 \Gamma^\infty$ with $V_0 (\Gamma^{\NN\mathbf{1}})^\infty$, and these identifications are compatible with the shift maps.

If $x \sim_{\Gg_{\Gamma^{\NN\mathbf{1}}}} y$, then $\sigma^p(x) \sim_{\Gg_{\Gamma^{\NN\mathbf{1}}}} \sigma^p(y)$ and similarly for $n$, so $\sigma^p$ and $\sigma^n$ descend to continuous maps $\tilde\sigma^p : q(V_0\Gamma^{\NN\mathbf{1}})^\infty \to q(V_p\Gamma^{\NN\mathbf{1}})^\infty$ and $\tilde\sigma^n : q(V_p\Gamma^{\NN\mathbf{1}})^\infty \to q(V_0\Gamma^{\NN\mathbf{1}})^\infty$. Since $x \sim_{\Gg_{\Gamma^{\NN\mathbf{1}}}} \sigma^{p+n}(x) = \sigma^p(\sigma^n(x))$ we see that $\tilde\sigma^p \circ \tilde\sigma^n$ is the identity map on $q(V_p\Gamma^{\NN\mathbf{1}})^\infty$, and similarly $\tilde\sigma^p \circ \tilde\sigma^n$ is the identity map on $q(V_0\Gamma^{\NN\mathbf{1}})^\infty$. So $\tilde\sigma^p$ and $\tilde\sigma^n$ are mutually inverse, and hence homeomorphisms.
\end{proof}

\begin{cor}\label{cor:both sc->Hausdorff}
Let $\Gamma$ be a row-finite source-free $k$-graph. Suppose that both $\Gamma$ and the sub-$1$-graph $\Gamma^{\NN\mathbf{1}}$ are simply connected. Then $\Gamma^\infty/\Gg_\Gamma$ is Hausdorff.
\end{cor}
\begin{proof}
Proposition~\ref{thm:univisAF} gives $f: \Gamma^0 \to  \ZZ^k$ such that $d(\lambda) = f(s(\lambda)) -  f(r(\lambda))$ for all $\lambda \in \Gamma$. Let $E$ be the directed graph such that $E^0 = f^{-1}(\ZZ\mathbf{1})$ and $E^1 = E^0 \Gamma^{\mathbf{1}}$. Lemma~\ref{lem:orbit spaces homeomorphic} gives $\Gamma^\infty/\Gg_\Gamma \cong E^\infty/\Gg_E$. Since $E^*$ is a sub-1-graph of the simply connected graph $\Gamma^{\NN \mathbf{1}}$, it is simply connected. Hence $E^\infty/\Gg_E$ is Hausdorff by \cite[Lemma~4.2]{KP1} (see Remark~\ref{rmk:tree ends Hausdorff}) and thus $\Gamma^\infty/\Gg_\Gamma$ is Hausdorff.
\end{proof}

\begin{example}
Surprisingly, simple connectedness of $\Gamma$ and of $\Gamma^{\NN\mathbf{1}}$ are independent conditions. For the monoidal $2$-graph $\Lambda$ of \cite[Example 7.1]{PQR1} (Example \ref{ex:pqrex1}), we have an isomorphism $\Pi(\Lambda) \cong \ZZ^2$ that intertwines $i : \Lambda \to \Pi (\Lambda)$ with $d : \Lambda \to \NN^2 \subseteq \ZZ^2$. So $\Gamma := \ZZ^2 \times_d \Lambda \cong \pi ( \Lambda ) \times_i \Lambda$ is simply connected. But $\Gamma^{\NN\mathbf{1}}$ is the graph with vertices $\{v_m : m \in \ZZ^2\}$ and six parallel edges from $v_{m + \mathbf{1}}$ to $v_m$ for each $m \in \ZZ^2$, so is not simply connected. In the other direction, let $\Delta_1$ be the $1$-graph with vertices $\ZZ$ and edges $e_n$ with $s(e_n) = n+1$ and $r(e_n) = n$, and define $l : \NN^2 \to \NN$ by $l(m,n) = m+n$. Then the $2$-graph $\Gamma := l^*(\Omega_1)$ has fundamental group $\ZZ$ generated by $(e_0, (1,0))(e_0, (0,1))^{-1}$, so is not simply connected, but $\Gamma^{\NN\mathbf{1}}$ is a disjoint union of copies of $\Omega_1$, so is simply connected.
\end{example}

\begin{rmk}
In the context of Corollary~\ref{cor:both sc->Hausdorff}, simple connectedness of $\Gamma^{\NN\mathbf{1}}$ is equivalent to that of $E^*$ as in Lemma~\ref{lem:orbit spaces homeomorphic}. Also, as in the proof of Proposition~\ref{prp:Hausdorffness}, the orbit space $(\Gamma^{\NN\mathbf{1}})^\infty / \Gg_{\Gamma^{\NN\mathbf{1}}}$ is a topological disjoint union of copies of $E^\infty/\Gg_E$ indexed $\ZZ^k/\ZZ\mathbf{1}$.
\end{rmk}

\begin{proof}[Proof of Theorem~\ref{thm:mainC*result}(iii)]
As in the proof of (2), since $\Lambda \to \Pi(\Lambda)$ is injective, $\Sigma$ is singly connected, and $C^*(\Sigma)$ is type I$_0$. The proof of Proposition~\ref{prop:fell} shows that $\mathcal{G}_\Sigma$ has trivial isotropy. Hence the spectrum of $C^*(\Sigma)$ is homeomorphic to the orbit space $\Sigma^\infty/\Gg_\Sigma$ \cite[Corollary~4.2]{OrloffClark}. Now, since $\Sigma^{\NN\mathbf{1}}$ is simply connected, Corollary~\ref{cor:both sc->Hausdorff} implies that $\Sigma^\infty/\Gg_\Sigma$ is Hausdorff. So $C^*(\Sigma)$ is a continuous-trace $C^*$-algebra. Since $X := \Sigma^\infty/\Gg_\Sigma$ is zero-dimensional, $\check{H}_3(X, \mathbb{Z}) = \{0\}$, and hence the Dixmier--Douady invariant $\delta(C^*(\Sigma))) \in \check{H}_3(X, \mathbb{Z})$ is trivial. So by the Dixmier--Douady theorem \cite[Corollary~5.58]{tfb}, $C^*(\Sigma)$ is Rieffel--Morita equivalent to $C(\Sigma^\infty/\Gg_\Sigma)$.
\end{proof}

\begin{rmk}
A  related realisation of $C^*$-algebras of $k$-graphs (and more general categories) as crossed products of abelian algebras by \emph{partial} actions of their fundamental groups appears in \cite[Theorem~4.17]{BBD}. Interestingly, embeddability also crops up there, for different reasons.
\end{rmk}

\begin{rmk}\label{rmk:confused}
It seems hard to nail down the relationships between the key hypotheses in this section: simple connectedness of $\Gamma$ and of $\Gamma^{\NN\mathbf{1}}$, and embedding of $\Gamma$ in $\Pi(\Gamma)$.

For example the following two assertions both seem reasonable: that if $\Gamma$ is simply connected, then the $\mathbf{1}$-dual $\mathbf{1}\Gamma$ obtained from Proposition~\ref{prop:gpb} for $f : n \mapsto n + \mathbf{1}$ is also simply connected; and that $\mathbf{1}\Gamma$ always embeds in $\Pi(\mathbf{1}\Gamma)$ (after all, $\mathbf{1}\Gamma \owns \lambda \mapsto (r(\lambda), d(\lambda), s(\lambda))$ is injective on $\bigcup_{n \le \mathbf{1}} \Gamma^n$, and this map descends to $\Pi(\Gamma)$, so the skeleton and factorisation rules are preserved in $\Pi(\Gamma)$). But at most one of these assertions is true in general: consider the skew-product $\Gamma := \ZZ^2 \times_d \Lambda$ of Example~\ref{ex:pqrex1}; we show that if $\mathbf{1}\Gamma$ is simply connected, then it does not embed in $\Pi(\mathbf{1}\Gamma)$.

Since $\mathbf{1}\Gamma$ is canonically isomorphic to the skew-product $\ZZ^2 \times_d (\mathbf{1}\Lambda)$, if $\mathbf{1}\Gamma$ is simply connected, then $\ZZ^2 \times_d (\mathbf{1}\Lambda)$ is simply connected, forcing $\pi_1(\Lambda, v) \cong \ZZ^2$. But inspection of the skeleton of $\mathbf{1}\Lambda$ shows that $eeec$ and $eedec$ are distinct blue cycles based at the vertex $ec \in \mathbf{1}\Lambda$, so generate a sub-semigroup of $\mathbf{1}\Lambda$ isomorphic to $\mathbb{F}^+_2$, which cannot embed in $\ZZ^2$.
\end{rmk}

\begin{rmk}\label{rmk:more confused}
The preceding remark is exemplary of a number of seemingly elementary questions that we have been unable to resolve. For example:
\begin{enumerate}[(i)]
\item If $\Gamma$ is simply connected and embeds in $\Pi(\Gamma)$, is $\Gamma^\infty/\Gg_\Gamma$ Hausdorff?
\item If both $\Gamma$ and $\Gamma^{\NN\mathbf{1}}$ are simply connected, does $\Gamma$ necessarily embed in $\Pi(\Gamma)$?
\item Which, if either, of the two assertions mentioned in Remark~\ref{rmk:confused} is correct?
\item Does $\mathbf{1}\Gamma$ always embed in $\Pi(\mathbf{1}\Gamma)$?
\end{enumerate}
\end{rmk}

\section{\texorpdfstring{$\tilde{A_2}$}{A\_2-tilde}-groups} \label{sec:tg}

In this section we construct coverings $\Sigma_\Tt \to \Lambda_\Tt$ of $2$-graphs corresponding to $\tilde{A_2}$-groups $\Gamma_\Tt$. These groups arise from free, vertex-transitive actions on buildings. We show that $\Sigma_\Tt$ and $\Lambda_\Tt$ both embed in their fundamental groupoids, and that $\Sigma_\Tt$ is always singly connected so that its $C^*$-algebra is of Type I$_0$.

The $\tilde{A_2}$-groups are built from finite projective planes. A finite projective plane $(P,L)$ of order $q$ consists of finite sets $P$ of \emph{points} and $L$ of \emph{lines} with $|P| = |L| = q^2 + q + 1$, and a relation $\in$ from $P$ to $L$---if $p \in l$, we say $p$ \emph{lies on} $l$ and that $l$ \emph{contains} $p$---such that any two points lie on exactly one common line, any two lines contain exactly one common point, and there exist four distinct points of which no single line contains more than two. Each line necessarily contains exactly $q$ points and each point necessarily lies on exactly $q$ lines.

We begin with a brief introduction of the groups we wish to study, and by collecting some structural results that we will need for our construction.

\subsection{\texorpdfstring{$\tilde{A_2}$}{A\_2-tilde}-group basics} \label{sec:tg-intro}

Following \cite[\S 2]{CMSZ1} given a finite projective plane $(P,L)$ and a bijection $\lambda: P \to L$, we define a \emph{triella compatible with $\lambda$} to be a set $\Tt \subset P \times P \times P$ such that
\begin{enumerate}[(T1)]
\item given $x,y \in P$, there exists $z \in P$ such that $(x,y,z) \in \Tt$ if and only if $y \in \lambda (x)$;
\item $(x,y,z) \in \Tt \Rightarrow (y,z,x) \in \Tt$;
\item for any $x,y \in P$, there is at most one $z \in P$ such that $(x,y,z) \in \Tt$.
\end{enumerate}

\begin{dfn}
Given a finite projective plane  $(P,L)$, a bijection $\lambda: P \to L$, and a triella $\Tt$ compatible with $\lambda$  as above, we define the associated
\emph{$\tilde{A_2}$-group} by
\begin{equation} \label{eq:triangledef}
\Gamma = \Gamma_\mathcal{T} := \langle  a_x, x \in P \mid   a_x a_y a_z = 1 \text{ for each } (x,y,z) \in \mathcal{T}  \rangle .
\end{equation}
\end{dfn}

\begin{rmks}
\begin{enumerate}[(i)]
\item The associated $\tilde{A_2}$-building is an oriented simplicial $2$-complex constructed from the Cayley graph of $\Gamma_\Tt$: the vertices or 0-simplices are identified with $\Gamma_\Tt$, the $1$-simplices are identified with pairs $(w, wa_x)$ where $w \in \Gamma_\Tt$  and $x \in P$. The $2$-simplices are identified with triples $(w, wa_x, wa_xa_y)$ where $w \in \Gamma_\Tt$, $x \in P$ and $y \in \lambda (x)$. The free and transitive action of $\Gamma_\Tt$ on 0-simplices by left multiplication extends to a free action on the building.

\item In \cite{V, KV} Vdovina et al{.} start with similar data to produce an object they call a polyhedron satisfying rules that have the flavour of a triella. We discovered this point of view late in our investigation and plan to look deeper into it in future work.
\end{enumerate}
\end{rmks}

\begin{example}\label{eg:Gamma A.1}
Many examples are considered in \cite{CMSZ2}. The following illustrative example with $q = 2$ was first described in \cite[\S4]{CMSZ2}:
\[
\Gamma_{A.1} = \langle a_0 , \ldots , a_6: a_{[i]_7} a_{[i+1]_7} a_{[i+3]_7} =1 \rangle, \text{ where$[i]_7 = i \!\!\!\!\pmod{7}$.}
\]
\end{example}

We describe elements of $\Gamma_\Tt$ as products of generators and their inverses. The following standard terminology for finitely-generated groups helps us discuss such expressions.

\begin{dfn}
Let $\Gamma_\Tt$ be an  $\tilde{A_2}$-group with generators $\{ a_x : x \in P\}$. By a \emph{word} in $\Gamma_\Tt$ we mean a string of the form $g_1 g_2 \cdots g_k$ such that each $g_i \in \{a_x, a_x^{-1} : x \in P\}$. The word $g_1 \cdots g_k$ \emph{represents} the element $w \in \Gamma_\Tt$ if the product $\prod^k_{i=1} g_i$ in $\Gamma_\Tt$ is equal to $w$. We typically indicate the group law by juxtaposition, so we write $w = g_1 \cdots g_k$ when the word $g_1 \cdots g_k$ represents $w$. Context will dictate whether a string $g_1 \cdots g_k$ is being regarded as a word or as a product.
\end{dfn}

It is helpful to express elements of $\Gamma_\Tt$ in a standard form.

\begin{prop}\label{prop:wd}
Let $\Gamma_\mathcal{T}$ be an  $\tilde{A_2}$-group  with  generators $\{ a_x : x \in P \}$.  Let $w \in \Gamma_\mathcal{T}$.
Then there are unique integers $m, n \ge 0$ and unique elements $x_1 , \ldots , x_m, y_1 , \ldots, y_n \in P$ such that
\begin{equation} \label{eq:rnf}
w = a_{x_1} \cdots  a_{x_{m}} a_{y_1}^{-1} \cdots a_{y_{n}}^{-1},  \text{ and }
\end{equation}
\[
\text{(a) } x_{i+1} \not\in \lambda  (x_{i} ) \text{ for $1 \le i < m$;}\quad
\text{(b) } y_{j} \not\in  \lambda ( y_{j+1} ) \text{ for $1 \le j < n$; and}\quad
\text{(c) } x_m \neq y_1 \text{ if $m,n \geq 1$.}
\]
For the same $m,n$, there are also unique  $t_1 , \ldots , t_n, s_1 , \ldots, s_m \in P$ such that
\begin{equation} \label{eq:lnf}
w = a_{t_1}^{-1} \cdots  a_{t_{n}}^{-1} a_{s_1} \cdots a_{s_{m}},  \text{ and }
\end{equation}
\[
\text{(a$'$) } s_{i+1} \not\in \lambda (s_{i}) \text{ for $1 \le i < m$;}\quad
\text{(b$'$) } t_{j} \not\in \lambda ( t_{j+1} ) \text{ for $1 \le j < n$; and}\quad
\text{(c$'$) } t_{\ell} \neq s_1 \text{ if $m, n \geq 1$.}
\]
We call the expressions above the \emph{right normal form}  and \emph{left normal form} of $w$ respectively. Both have minimal length amongst words in the generators and their inverses that represent $w$.
Moreover, every minimal-length word in the generators and their inverses that represents $w$ contains $m$ generators and $n$ generator-inverses.
\end{prop}
\begin{proof}
See \cite[Proposition 3.2]{CMSZ1} and \cite[Lemma 6.2]{CM}.
\end{proof}

\begin{cor}\label{cor:diamond}
Let $\Gamma_\Tt$ be an  $\tilde{A_2}$-group  with generators $\{ a_x : x \in P\}$. For all $x , y \in P$ such that $x \neq y$,
there exist unique  $s, t \in P$ with $s \ne t$ such that $a_x^{-1}  a_y= a_s a_t^{-1}$.
\end{cor}

Proposition~\ref{prop:wd} allows us to define a degree functor for a $2$-graph structure on $\Gamma_\Tt$ in terms of the number of generators and their inverses in a minimal representative of an element.

\begin{dfn} \label{dfn:delta}
Let $\Gamma_\mathcal{T}$ be an  $\tilde{A_2}$-group . Define $\delta : \Gamma_\mathcal{T} \to \NN^2$ by $\delta (w) = (m,n)$ if
its right normal form is as in equation~\eqref{eq:rnf} (equivalently, its left normal form  is as in equation~\eqref{eq:lnf}).
We define the \emph{length} of $w$ to be $|\delta(w)| = m + n$. We call $\delta$ the \emph{shape function}.
\end{dfn}

\begin{rmk} \label{rmk:deltanothom}
The shape function $\delta$ is not additive: in $\Gamma_{A.1} = \langle a_0 , \ldots , a_6: a_{[i]_7} a_{[i+1]_7} a_{[i+3]_7} =1 \rangle$,
\[
 \delta (a_1 a_2) = \delta ( a_4^{-1} ) = (0,1) \neq (2,0) = \delta ( a_1 ) + \delta ( a_2 ).
\]
\end{rmk}

The shape function $\delta$ gives rise to a natural notion of a reduced word.

\begin{dfn}
A word $g_1 \cdots g_k$ in $\Gamma_{\mathcal{T}}$ is said to be
\emph{reduced} if it has minimal length among words that represent the same element of $\Gamma_\Tt$. That is,
$g_1 \cdots g_k$ is reduced if $|\delta(g_1 \cdots g_k)| = k$.
\end{dfn}

\begin{rmks}\label{rmks:normal vs reduced}
\begin{enumerate}[(i)]
\item The final statement of Proposition~\ref{prop:wd} shows that words in right normal form or left normal form are reduced words.
\item Not all words that have no ``obvious cancellations" are reduced: the word $w_1 = a_0{a_4}^{-1}a_6$ in $\Gamma_{A.1} = \langle a_0 , \ldots , a_6: a_{[i]_7} a_{[i+1]_7} a_{[i+3]_7} =1 \rangle$,
is not reduced since
    \[
        a_0 {a_4}^{-1} a_6 =   a_0 a_1a _2 a_6 = a_0 a_1{a_0}^{-1} = {a_3}^{-1}{a_0}^{-1}.
    \]
\item Every subword of a reduced word is reduced.
\item If  $w = g_1 \cdots g_k$ is reduced and for some $i$, $g_i = a_x$ and $g_{i+1}= a_y^{-1}$ for some $x, y \in P$ with $x \ne y$, then by Corollary~\ref{cor:diamond}, there exist unique $s, t \in P$ with $s \ne t$ such that $a_x{a_y}^{-1} = {a_s}^{-1}a_t$. The word obtained from $w$ by replacing $g_ig_{i+1}= a_x a_y^{-1}$  with ${a_s}^{-1}a_t$ is also reduced.
\end{enumerate}
\end{rmks}

\begin{example}
Consider $\Gamma_\Tt := \Gamma_{A.1} = \langle a_0 , \ldots , a_6: a_{[i]_7} a_{[i+1]_7} a_{[i+3]_7} =1 \rangle$ from Example~\ref{eg:Gamma A.1}.
For $w = a_0 a_2^{-1} a_5^{-1}  \in \Gamma_\Tt$, we have $\delta(w) = (1,2)$; the reduced expressions for $w$, and the corresponding segment of the reversed Cayley
graph of $\Gamma_\Tt$ (the Cayley graph of $\Gamma_\Tt^{op}$), are illustrated below.
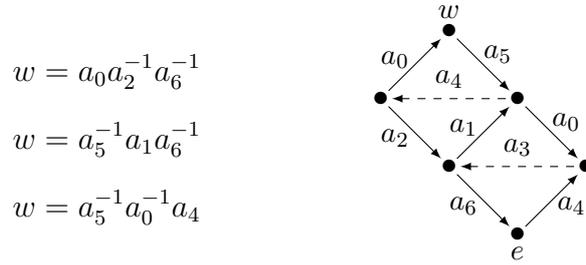
\[
\begin{tikzpicture}[xscale=.9,yscale=.9]

\node at (-6,2.4) {$w=a_0 a_2^{-1} a_6^{-1}$};
\node at (-6,1.4) {$w=a_5^{-1} a_1 a_6^{-1} $};
\node at (-6,0.4) {$w=a_5^{-1} a_0^{-1} a_4$};

\node[circle, inner sep=0pt] (v00) at (0,0) {$\bullet$};

\node[circle, inner sep=0pt] (v01) at (-1,1) {$\bullet$};
\node[circle, inner sep=0pt] (v10) at (1,1) {$\bullet$};

\node[circle, inner sep=0pt] (v11) at (0,2) {$\bullet$};
\node[circle, inner sep=0pt] (v02) at (-2,2) {$\bullet$};

\node[circle, inner sep=0pt] (v12) at (-1,3) {$\bullet$};

\node[anchor=south, inner sep=0pt] at (v12.north) {\smash{\small $w$}};
\node[anchor=north, inner sep=5pt] at (v00.south) {\smash{\small $e$}};

\draw[->,-latex] (v00) -- (v10) node[pos=0.4,right] {\small $a_4$};
\draw[->,-latex] (v01) -- (v00) node[pos=0.6,left] {\small $a_6$};

\draw[->,-latex] (v01) -- (v11) node[pos=0.6,left] {\small $a_1$};
\draw[->,-latex] (v11) -- (v10) node[pos=0.3,right] {\small $a_0$};

\draw[->,-latex] (v02) -- (v01) node[pos=0.6,left] {\small $a_2$};
\draw[->,-latex] (v02) -- (v12) node[pos=0.6,left] {\small $a_0$};

\draw[->,-latex] (v12) -- (v11) node[pos=0.3,right] {\small $a_5$};
\draw[latex-,dashed] (v01) -- (v10) node[pos=0.5,above] {\small $a_3$};
\draw[latex-,dashed] (v02) -- (v11)node[pos=0.5,above] {\small $a_4$};

\end{tikzpicture}
\]
\end{example}

To obtain $2$-graphs from $\tilde{A_2}$ groups, we relate the shape function $\delta$ to the group law.

\begin{lem}[Unique factorisation] \label{lem:deltadd} 
Let $\Gamma_\mathcal{T}$ be an  $\tilde{A_2}$-group  and suppose that $m,n \in \NN^2$ and $w \in \Gamma_\mathcal{T}$ satisfy $\delta (w) = m+n$. Then there exist unique $h,k \in \Gamma_\mathcal{T}$ such that  $\delta (h)= m$, $\delta (k)=n$ and $w=hk$.
More generally, if $n_i \in \NN^2$ satisfy $\delta (w) = n_1 + \cdots + n_k$, then there exist unique $h_i \in \Gamma_\mathcal{T}$ such that each $\delta (h_i)= n_i$ and $w= h_1 \cdots h_k$.

Given $w, h, k  \in \Gamma_\mathcal{T}$ such that $\delta(whk) = \delta(w) + \delta(h)  + \delta(k)$,
we have $\delta(wh) = \delta(w) + \delta(h)$ and $\delta(hk) = \delta(h) + \delta(k)$.
\end{lem}

\begin{proof}
This follows from repeated applications of Corollary~\ref{cor:diamond}.
\end{proof}

\begin{nota} \label{rmk:end}
If $\delta (w) = (m,n) \ge \mathbf{1}$ then Lemma~\ref{lem:deltadd} yields unique $a,b,c,d\in \Gamma_\mathcal{T}$ such that
\begin{equation} \label{eq:topandbottom}
w = bd = c a, \quad \delta ( a ) = \delta( b ) = \mathbf{1} \quad\text{and} \quad \delta ( d ) = \delta (c) = \delta (w) - \mathbf{1}.
\end{equation}
\noindent
We adopt the notation $s(w)=a$,  $r(w)=b$, $c(w)=c$, $d(w)=d$.
Note that if  $\delta (w) = \mathbf{1}$, then $r (w)= s (w)$ and $b(w) = c(w)= 1$.
\end{nota}

\noindent
We provide a criterion for determining when a concatenation of three reduced words is reduced.

\begin{prop}\label{prop:reduced-criterion}
Let $\Gamma_\mathcal{T}$ be an  $\tilde{A_2}$-group  and fix $w_0, w_1, w_2 \in \Gamma_\mathcal{T}$.  Suppose that $\delta(w_0w_1) = \delta(w_0) + \delta(w_1)$,
$\delta(w_1w_2) = \delta(w_1) + \delta(w_2)$ and $\delta(w_1) \ge \mathbf{1}$.
Then
\[
\delta(w_0w_1w_2) = \delta(w_0) + \delta(w_1) + \delta(w_2).
\]
\end{prop}

\begin{proof}
We induct on $|\delta(w_2)|$. Suppose that $|\delta(w_2)|=1$, so $\delta(w_2) \in \{(1,0), (0,1)\}$.
If $\delta(w_2) = (1, 0)$, then $w_2 = a_x$ for some $x \in P$, so $w_0w_1w_2 = w_0w_1a_x$.
By Proposition~\ref{prop:wd}, if $\delta( w_0w_1 ) = (m, n)$, then in left normal form,
$w_0w_1 = a_{s_1}^{-1} \cdots  a_{s_{n}}^{-1} a_{t_1} \cdots a_{t_{m}}$
and $w_1 = a_{p_1}^{-1} \cdots  a_{p_{k}}^{-1} a_{q_1} \cdots a_{q_{\ell}}$.
Lemma~\ref{lem:deltadd} gives $q_{\ell} = t_m$ (as $\delta(w_0w_1) = \delta(w_0) + \delta(w_1)$)
and $x \notin \lambda(q_{\ell})$ (as $\delta(w_1w_2) = \delta(w_1) + \delta(w_2)$).
Hence
\[
w_0w_1w_2= a_{s_1}^{-1} \cdots  a_{s_{n}}^{-1} a_{t_1} \cdots a_{t_{m}}a_x
\]
is the left normal form of $w_0w_1w_2$ and so
\[
\delta(w_0w_1w_2) = (m+1, n)= (m, n) + (1, 0) = \delta(w_0w_1) + \delta(w_2)= \delta(w_0) + \delta(w_1) + \delta(w_2).
\]
If $\delta(w_2) = (0,1)$, arguing similarly with right normal forms gives
$\delta(w_0w_1w_2) = \delta(w_0) + \delta(w_1) + \delta(w_2)$.

Now suppose that the result holds for $|\delta(w_2)| = n \ge 1$, suppose that $|\delta(w_2)| = n + 1$.
Then there exist unique $h, k  \in \Gamma_\mathcal{T}$ such that $w_2 = hk$, $\delta(w_2) = \delta(h) + \delta(k)$ and $|\delta(k)|= 1$.
Since
\[
\delta(w_1hk) = \delta(w_1w_2) = \delta(w_1) + \delta(w_2) = \delta(w_1) + \delta(h) + \delta(k),
\]
Lemma~\ref{lem:deltadd} gives $\delta(w_1h) = \delta(w_1) + \delta(h)$.
Since $|\delta(h)| = n$, the induction hypothesis gives
$\delta(w_0(w_1h))  = \delta(w_0) + \delta(w_1) + \delta(h) = \delta(w_0) + \delta(w_1h)$.
Moreover,
\begin{align*}
 \delta((w_1h)k) &= \delta(w_1w_2) = \delta(w_1) + \delta(w_2) = \delta(w_1) + \delta(h) + \delta(k) = \delta(w_1h) + \delta(k).
 \intertext{Therefore, since $\delta(w_1h) \ge \mathbf{1}$ and $|\delta(k)| = 1$, it follows that}
 \delta(w_0w_1w_2) &=  \delta(w_0(w_1h)k) =  \delta(w_0) + \delta(w_1h) + \delta(k) =  \delta(w_0) + \delta(w_1) + \delta(h)+ \delta(k) \\
      &= \delta(w_0) + \delta(w_1) + \delta(w_2). \qedhere
\end{align*}
\end{proof}

\noindent
The following extends the above criterion to an arbitrary concatenation of reduced words.

\begin{cor}\label{prop:reduced-tower}
 Fix $w_0, w_1, \dots, w_n \in \Gamma_\mathcal{T}$. Suppose that  $\delta(w_iw_{i+1}) = \delta(w_i) + \delta(w_{i+1})$ for all $0 \le i < n$ and $\delta(w_i) \ge \mathbf{1}$ for all $0 < i < n$. Then
\[
\delta(w_0w_1 \cdots w_{n-1} w_n) = \delta(w_0) + \delta(w_1) + \cdots + \delta(w_{n-1}) + \delta(w_{n}).
\]
\end{cor}

\begin{proof}
We induct on $n$. This is trivial for $n = 1$. Fix $n \ge 1$, suppose the result holds for all $k \le n$, and fix
$w_0, w_1, \dots, w_n, w_{n+1} \in \Gamma_\mathcal{T}$ with $\delta(w_iw_{i+1}) = \delta(w_i) + \delta(w_{i+1})$ for all $i = 0, 1, \dots, n$ and  $\delta(w_i) \ge \mathbf{1}$ for all $i = 1, \dots, n$.
Then
\[
\delta(w_0w_1 \cdots w_n) =  \delta(w_0) + \delta(w_1) + \cdots + \delta(w_n) =  \delta(w_0w_1 \cdots w_{n-1})  + \delta(w_{n}),
\]
and since $\delta(w_n w_{n+1}) = \delta(w_n)+\delta(w_{n+1}) $, we have
\[
\delta(w_0w_1 \cdots w_nw_{n+1}) = \delta(w_0w_1 \cdots w_{n-1}) +  \delta(w_{n})  + \delta ( w_{n+1}) = \delta(w_0) + \delta(w_1) + \cdots + \delta(w_n) + \delta ( w_{n+1})
\]
by Proposition~\ref{prop:reduced-criterion}.  Thus the result holds by induction.
\end{proof}

\subsection{The \texorpdfstring{$2$}{2}-graph associated to an  \texorpdfstring{$\tilde{A_2}$}{A\_2-tilde}-group } \label{sec:2g}

Given an  $\tilde{A_2}$-group  $\Gamma_\mathcal{T}$, we now construct a $2$-graph $\Lambda_\mathcal{T}$ using the relation between the  multiplicative structure of its reduced words and the shape function discussed in the previous section.

\begin{dfn}\label{dfn:LambdaT}
Fix an  $\tilde{A_2}$-group  $\Gamma_\mathcal{T}$. We define
\[
\Lambda_\mathcal{T} = \{ w \in \Gamma_\mathcal{T} : \delta (w) \ge \mathbf{1} \}  \quad\text{and } \quad
\Lambda_\mathcal{T}^0 =  \{ u \in \Gamma_\mathcal{T} : \delta (u) = \mathbf{1} \} .
\]
We define $r, s : \Lambda_\Tt \to \Lambda_\Tt^0$ as in Notation~\ref{rmk:end} and $d : \Lambda_\mathcal{T} \to \NN^2$  by $d(\lambda) = \delta(\lambda) -\mathbf{1}$. For $\lambda,\mu \in \Lambda_\mathcal{T}$ such that $s(\lambda) = r(\mu)$,  we define $\lambda \circ \mu$ as follows: Write $\lambda = c(\lambda) s(\lambda)$ and $\mu = r (\mu) b(\mu)$ as in Notation~\ref{rmk:end}; we define
\begin{equation} \label{eq:circdef}
\lambda \circ \mu :=  c(\lambda)s (\lambda) b(\mu) .
\end{equation}
\end{dfn}

Our definition of $\lambda \circ \mu$ in~\eqref{eq:circdef} emphasises the overlap of $\lambda = c(\lambda)s(\lambda)$ and $\mu = r(\mu)c(\mu)$ in the element $s(\lambda) = r(\mu)$ of $\delta^{-1}(\mathbf{1}) \subseteq \Gamma_\Tt$. We can also express it to emphasise its compatibility with the maps $b$ and $c$: for $\lambda, \mu \in \Lambda_\Tt$ as above with $s(\lambda) = r(\mu)$,
\begin{equation} \label{eq:circequiv}
\lambda \circ \mu = c(\lambda)s(\lambda)b(\mu) = r(\lambda) b(\lambda) b(\mu)  \quad\text{ and }\quad
\lambda \circ \mu = c(\lambda) r ( \mu ) b(\mu) = c(\lambda) c(\mu) s(\mu).
\end{equation}

Our main result in this subsection is that Definition~\ref{dfn:LambdaT} defines a $2$-graph:

\begin{thm} \label{thm:lambdaT}
With  definitions and notation as above, $(\Lambda_\Tt, d )$  is a $2$-graph, and the maps $b, c : \Lambda_\mathcal{T} \to \Gamma_\Tt$ of Notation~\ref{rmk:end} are $1$-cocycles.
\end{thm}

\begin{proof}
Associativity of multiplication in $\Gamma_\Tt$ ensures that $\Lambda_\Tt$ is a category under $\circ$. To see that $d : \Lambda_\Tt \to \NN^2$ is a functor,
fix $\lambda,\mu \in \Lambda_\Tt$ with $s(\lambda)=r(\mu)$. We have $\lambda \circ \mu= c(\lambda)s(\lambda)b(\mu)$, where $\delta(s(\lambda))) = \mathbf{1}$.
So the first part of Lemma~\ref{lem:deltadd} gives
\begin{align*}
\delta(c(\lambda)s(\lambda)) &= d(\lambda) + \mathbf{1} = \delta(c(\lambda)) + \delta(s(\lambda)) \\
\delta(s(\lambda)b(\mu)) &=  d(\mu) +  \mathbf{1} = \delta(s(\lambda)) + \delta(b(\mu)).
\end{align*}
Hence by Proposition~\ref{prop:reduced-criterion}, and since $\delta(s(\lambda)) = \mathbf{1}$ by definition,
\[
d(\lambda \circ \mu) = \delta(c(\lambda)s(\lambda) b(\mu)) - \mathbf{1}
    =   \delta(c(\lambda)) + \delta(s(\lambda)) + \delta( b(\mu)) - \mathbf{1}
    = d(\lambda) + d(\mu).
\]

It remains to show that $( \Lambda_\mathcal{T},d)$ satisfies the factorisation property. Suppose that $d(\lambda) = (m_1+m_2,n_1+n_2)$.
Then $\delta ( \lambda ) = (m_1+m_2,n_1+n_2)+\mathbf{1}$. Hence by Lemma~\ref{lem:deltadd} there exist unique $g, h, k$  such that $ \lambda = ghk$, $\delta(g) = (m_1, n_1)$, $\delta(h) = \mathbf{1}$ and $\delta(k) = (m_2, n_2)$.  Thus, $\lambda = \mu \circ \nu$ where $\mu = gh$, $\nu = hk$, $d(\mu) = (m_1, n_1)$, and $d(\nu) = (m_2, n_2)$ and this is the unique such factorisation.

Fix $\lambda, \mu \in \Lambda_\mathcal{T}$ with $s(\lambda) = r(\mu)$. Equation~\ref{eq:circequiv} and the definition of $c$ give $c(\lambda \circ \mu)s(\lambda\circ\mu) = \lambda\circ\mu = c(\lambda)c(\mu)s(\mu)$ and $r(\lambda\circ\mu)b(\lambda\circ\mu) = \lambda\circ\mu = r(\lambda) b(\lambda)b(\mu)$. We already saw that $s(\lambda\circ\mu) = s(\mu)$ and $r(\lambda\circ\mu) = r(\lambda)$, so cancellativity in $\Gamma_\Tt$ gives $c(\lambda)c(\mu) = c(\lambda\circ\mu)$ and $b(\lambda)b(\mu) = b(\lambda\circ\mu)$.
\end{proof}

\begin{rmk}\label{rmk:us vs RobSteg}
Resume the notation of \cite[\S7]{RobSteg1}. Let $M_1, M_2$ be the matrices \cite[p.135]{RobSteg1} obtained from the Cayley graph $\mathscr{B}_\Tt$ of $\Gamma_\Tt$ regarded as a building as in \cite{CMSZ1}. Then
$\Lambda_\Tt$ is isomorphic to the $2$-graph $\Lambda_{M_1, M_2}$ obtained from the $M_i$ as in \cite[Examples~1.7(iv)]{KP2}. Indeed, as $\Gamma_\Tt$ acts transitively on vertices of $\mathscr{B}_\Tt$, we can identify the alphabet $A = \Gamma/\mathcal{I}$ \cite[p.135]{RobSteg1} with type-rotating isometries $i : t \to \mathscr{B}_\Tt$ such that $i((0,0)) = e_{\Gamma_\Tt}$. By Proposition~\ref{prop:wd}, $i \mapsto i((1,1))$ is a bijection between such isometries and $\delta^{-1}(\mathbf{1}) = \Lambda_\Tt^0$. Likewise, for $w_1, w_2 \in \Lambda_\Tt^0$, the set $w_1\Lambda_\Tt^{\varepsilon_i}w_2$ is in bijection with type-rotating isometries $i : \mathfrak{p}_{\varepsilon_i} \to \mathscr{B}_\Tt$ such that $i((0,0)) = e_{\Gamma_\Tt}$, $t(\mathbf{1}) = w_1$, and $t(\mathbf{1}+\varepsilon_i)t(\varepsilon_i)^{-1} = w_2$; that is, diagrams as in \cite[Figure~9]{RobSteg1}. So the adjacency matrices of $\Lambda_\Tt$ are the $M_i$. Since $M_1, M_2$ satisfy (H0)--(H3) \cite[Proposition~7.9 and Theorem~7.10]{RobSteg1}, $M_1M_2$ is a $0,1$-matrix, so \cite[Theorems 4.4~and~4.5]{HRSW} gives $\Lambda_\Tt \cong \Lambda_{M_1, M_2}$.
\end{rmk}

\begin{cor}\label{cor:LambdaT embeds}
With notation as above, the cocycle $c: \Lambda_\mathcal{T} \to \Gamma_\mathcal{T}$ of Theorem~\ref{thm:lambdaT} is essential and the canonical map
$i : \Lambda_\mathcal{T} \to \Pi(\Lambda_\mathcal{T})$ is injective.
\end{cor}
\begin{proof}
Since $\lambda = c(\lambda) s(\lambda)$ for all $\lambda \in \Lambda_\mathcal{T}$,
$c \times s$ is injective. Hence $c$ is essential as in Definition~\ref{def:essential}, and the result follows from Proposition~\ref{lprop:ess-emb}.
\end{proof}

\subsection{The covering \texorpdfstring{$2$}{2}-graph\texorpdfstring{ $\Sigma_\Tt$}{}}\label{sec:cover-2g}
In this section we construct a covering $2$-graph $\Sigma_\Tt$ for $\Lambda_\Tt$.

We define $\Sigma_\Tt \subseteq \Gamma_\Tt \times \Gamma_\Tt$ as follows. Let
\[
\Sigma_\Tt := \{ (x, y) \in \Gamma_\Tt \times \Gamma_\Tt : \mathbf{1} \le \delta(x^{-1}y)  \} \quad\text{and}  \quad
\Sigma^0_\Tt := \{ (x, y) \in \Gamma_\Tt \times \Gamma_\Tt : \mathbf{1} = \delta(x^{-1}y)   \}
\]

\noindent
with  $d(x, y)  := \delta(x^{-1}y) - \mathbf{1}$ for all $(x, y) \in \Sigma_\Tt$.
By  Lemma~\ref{lem:deltadd} for each $(x, y) \in  \Sigma_\Tt$ there exist unique $z_{x,y}, w_{x,y} \in \Gamma_\Tt$ such that
\begin{equation}\label{eq:r and s}
\delta(x^{-1}z_{x,y}) = \delta(w_{x,y}^{-1}y) = \mathbf{1} \quad\text{and} \quad \delta(x^{-1}y) = \delta(x^{-1}z_{x,y}) + \delta(z_{x,y}^{-1}y) = \delta(x^{-1}w_{x,y}) + \delta(w_{x,y}^{-1}y).
\end{equation}
We define $r(x, y) := (x, z_{x,y})$ and $s(x, y) := (w_{x,y}, y)$. If $(u, v) \in \Sigma_\Tt$ satisfies $s(x, y) = r(u, v)$, we define $(x, y)(u, v) := (x, v)$.
We show that $\Sigma_\Tt \cong \Gamma_\Tt \times_{\tilde{c}} \Lambda_\Tt$ (see Definition~\ref{eg:skewproducts}).

\begin{prop}\label{prop:skew-cover}
With the above structure $\Sigma_\Tt$ is a $2$-graph. Let $c : \Lambda_\Tt \to \Gamma_\Tt$ be the cocycle of Theorem~\ref{thm:lambdaT}. There is an isomorphism  $\phi : \Sigma_\Tt \to \Gamma_\Tt \times_{c} \Lambda_\Tt$ such that $\phi (x, y) = (x, x^{-1}y)$ for all $(x, y) \in \Sigma_\Tt$. The inverse satisfies $\phi^{-1}(x, \lambda) = (x, x\lambda)$. There is a free action of $\Gamma_\Tt$ on $\Sigma_\Tt$ given by $g \cdot (x, y) := (gx, gy)$, and $\phi$ is equivariant for this action and the left action of $\Gamma_\Tt$ on $\Gamma_\Tt \times_{c} \Lambda_\Tt$ by translation in the first coordinate. In particular, $\phi$ descends to an isomorphism $\widetilde{\phi} : \Gamma_\Tt \backslash \Sigma_\Tt \to \Lambda_\Tt$ such that $\widetilde{\phi}([x,y]) = x^{-1}y$.
\end{prop}
\begin{proof}
We prove that $\psi:  \Gamma_\Tt \times_{c} \Lambda_\Tt \to \Sigma_\Tt$ defined by $\psi(x, \lambda) =  (x, x\lambda)$ for $(x, \lambda)  \in \Gamma_\Tt \times_{c} \Lambda_\Tt$ is an isomorphism, and deduce that $\phi = \psi^{-1}$ has the desired properties. To see that $\psi$ is a functor, fix $(x, \lambda) \in  \Gamma_\Tt \times_{c} \Lambda_\Tt$. Then  $d(x, \lambda) = d(\lambda) = d(x, x\lambda)= d(\psi(x, \lambda) )$. Recall that $c(\lambda)s(\lambda) = \lambda = r(\lambda)b(\lambda)$, $s(x, \lambda) = (xc(\lambda), s(\lambda))$ and $s(x, x\lambda) = (xc(\lambda), x\lambda)$ since
$w = xc(\lambda)$ satisfies $\delta(w^{-1}x\lambda) = \delta(s(\lambda)) = \mathbf{1}$ and
\[
\delta(x^{-1}x\lambda) = \delta(\lambda) =\delta(c(\lambda)) +  \delta(s(\lambda)) = \delta(x^{-1}w) + \delta(w^{-1}x\lambda).
\]
Hence
\[
\psi(s(x, \lambda)) = \psi(xc(\lambda), s(\lambda)) = (xc(\lambda), xc(\lambda)s(\lambda)) = (xc(\lambda), x\lambda) = s(x, x\lambda) = s(\psi(x, \lambda)).
\]
A similar computation shows that $\psi(r(x, \lambda)) = (x, xr(\lambda)) = r(\psi(x, \lambda))$.

Given composable elements $(x, \lambda), (xc(\lambda), \mu) \in \Gamma_\Tt \times_{c} \Lambda_\Tt$, the above argument shows that $\psi(x, \lambda), \psi(xc(\lambda), \mu)$ are composable in $\Sigma_\Tt$ .  We have
\begin{align*}
\psi((x, \lambda) (xc(\lambda), \mu))
    &=  \psi(x, \lambda\circ\mu)
     =  (x, x(\lambda\circ\mu))\\
    &= (x, xc(\lambda)\mu)
     = (x, x\lambda)(xc(\lambda), xc(\lambda)\mu)
     = \psi(x, \lambda) \psi(xc(\lambda), \mu).
\end{align*}
Hence, $\psi$ is an isomorphism and thus $\Sigma_\Tt$ is a $2$-graph. That $\phi$ is equivariant follows from its definition and the last assertion follows from  \cite[Remark~5.6]{KP2}.
\end{proof}

\begin{prop}  \label{prop:cover-1cnc}
Let $(x, z), (w, y) \in \Sigma_\Tt^0$.  Then $(x, z)\Sigma_\Tt(w, y) \not= \emptyset$ if and only if
\[
\delta(x^{-1}w) + \mathbf{1}  = \delta(x^{-1}y) =   \delta(z^{-1}y)  + \mathbf{1},
\]
and then $(x, z)\Sigma_\Tt(w, y) = \{(x, y)\}$. In particular, $\Sigma_\Tt$ is singly connected and $C^*(\Sigma_\Tt)$ is type~I$_0$.
\end{prop}
\begin{proof}
If $\sigma \in (x, z)\Sigma_\Tt(w, y)$, then $r(\sigma) = (x, z)$ and $s(\sigma) = (w, y)$,
so $\sigma = (x, y)$, $w = w_{x, y}$ and $z = z_{x, y}$ by~\eqref{eq:r and s}.
In particular, $(x, z)\Sigma_\Tt(w, y)$ is either empty or equal to $\{(x, y)\}$.

If $\delta(x^{-1}y) \ge \mathbf{1}$, then $(x,y) \in \Sigma_\Tt$ if and only if $s(x, y) = (w, y)$ and $r(x, y) = (x, z)$.
Moreover, $s(x, y) = (w, y)$ if and only if $w = w_{x, y}$, that is (see \eqref{eq:r and s})
\begin{align*}
\delta(x^{-1}y)  &= \delta(x^{-1}w) + \delta(w^{-1}y)  = \delta(x^{-1}w) + \mathbf{1}
\intertext{and  $r(x, y) = (x, z)$ if and only if  $z = z_{x, y}$, that is }
\delta(x^{-1}y)  &=   \delta(x^{-1}z)  + \delta(z^{-1}y)  =   \delta(z^{-1}y)  + \mathbf{1}.
\end{align*}
The final assertion follows from the first paragraph of the proof and Proposition~\ref{prop:fell}.
\end{proof}

\begin{rmk}
 That $\Sigma_\Tt$ is singly connected also follows from the facts that  $\Sigma_\Tt \cong \Gamma_\Tt \times_{\tilde{c}} \Lambda_\Tt$
 (by Proposition~\ref{prop:skew-cover}), $c$ is essential and  $\Gamma_\Tt \times_{\tilde{c}} \Lambda_\Tt$ is singly connected
 (by Corollary~\ref{cor:LambdaT embeds}).
\end{rmk}


\begin{thebibliography}{APCaHR}
\bibitem
{APCaHR} G. Aranda Pino, J. Clark, A. an Huef and I. Raeburn, \emph{Kumjian--Pask algebras of higher-rank graphs}, Trans. Amer. Math. Soc., \textbf{365} (2013), 3613--3641.

\bibitem
{AB} B. Armstrong and N. Brownlowe, \emph{Product-system models for twisted $C^*$-algebras of topological higher-rank graphs}, J. Math. Anal. Appl. \textbf{466} (2018),  1443--1475.

\bibitem
{APS} S. Allen, D. Pask and A. Sims, \emph{A dual graph construction for higher-rank graphs, and $K$-theory for finite $2$-graphs}, Proc. Amer. Math. Soc. \textbf{ 134} (2006), 455--464.

\bibitem
{BPRS} T.  Bates, D.  Pask, I.  Raeburn and W. Szyma\'{n}ski, \emph{The $C^*$-algebras of row-finite graphs}, New York J. Math.  \textbf{ 6} (2000), 307--324.

\bibitem{BKQ21} E. B\'edos, S. Kaliszewski, and J. Quigg, \emph{Skew products of finitely aligned left cancellative small categories and Cuntz--Krieger algebras}, M\"unster J. Math. \textbf{14} (2021), 59--99.

\bibitem
{BH} M. Bridson and A. Haefliger. \emph{Metric spaces of non-positive curvature}. 1999, Springer-Verlag ,Berlin.

\bibitem{BBD} K.A. Brix, C. Bruce and A. Dor On, \emph{Normal coactions extend to the $C^*$-envelope}, preprint 2023 (arXiv:2309.04817 [math.OA]).

\bibitem
{CM} D. Cartwright and W. M\l otkowski, \emph{Harmonic analysis for groups acting on triangle buildings}, J. Austral. Math. Soc. (Series A) \textbf{56} (1994), 345--383.


\bibitem
{CMSZ1} D. Cartwright, A. Mantero, T. Steger and A. Zappa, \emph{Groups acting simply transitively on the vertices of a building of type $\widetilde{A}_2$ I}, Geometric\ae Dedicata \textbf{47} (1993), 143--166.

\bibitem
{CMSZ2} D. Cartwright, A. Mantero, T. Steger and A. Zappa, \emph{Groups acting simply transitively on the vertices of a building of type $\widetilde{A}_2$ II}, Geometric\ae Dedicata \textbf{47} (1993), 167--226.

\bibitem
{OrloffClark} L.O. Clark, \emph{Classifying the types of principal groupoid $C^*$-algebras}, J. Operator Th. \textbf{57} (2007), 251--266.

\bibitem
{CaHS} L.O. Clark, A. an Huef and A. Sims, \emph{AF-embeddability of $2$-graph algebras and quasidiagonality of {$k$}-graph algebras}, J. Funct. Anal. \textbf{271}, (2016), 958--991.

\bibitem
{CFaH} L. Orloff Clark,  C. Flynn and A. an Huef, \emph{Kumjian-Pask algebras of locally convex higher rank graphs}, J. Algebra \textbf{399} (2014), 445--474.


\bibitem
{DPY}  K.R. Davidson, S.C. Power and D. Yang, \emph{Atomic representations of rank 2 graph algebras},  J. Funct. Anal. \textbf{255} (2008) 819--853.

\bibitem{Dehornoy} P. Dehornoy with F. Digne, E. Godelle, D. Krammer, and J. Michel, Foundations of Garside theory. EMS Tracts in Mathematics, 22. European Mathematical Society (EMS), Z\"urich, 2015.

\bibitem
{E} D. G. Evans, {\em On the $K$-theory of higher rank graph $C^*$-algebras}, New York J. Math. \textbf{14} (2008), 1--31.

\bibitem
{FPS} C. Farthing, D. Pask and A. Sims, \emph{Crossed products of $k$-graph algebras by $\mathbb{ Z}^\ell$}, Houston J. Math. \textbf{ 35} (2009), 903--933.

\bibitem
{G}  E. Gillaspy, \emph{K-theory and homotopies of $2$-cocycles on higher-rank graphs}, Pacific J. Math. \textbf{278} (2015),  407--426.

\bibitem
{HRSW} R. Hazlewood, I. Raeburn, A. Sims and S.B.G. Webster, \emph{On some fundamental results about higher-rank graphs and their $C^*$-algebras}, Proc. Edinburgh Math. Soc. \textbf{56} (2013), 575--597.

\bibitem{Higgins} P. J. Higgins, Notes on categories and groupoids, Van Nostrand Rienhold Mathematical Studies, no. 32, Van Nostrand Reinhold, London--New Yorkp--Melbourne, 1971.

\bibitem
{J} P.T. Johnstone, \emph{On embedding categories in groupoids}, Math. Proc. Cambridge Philos. Soc., \textbf{145}  (2008), 273--294.

\bibitem
{KKQS} S. Kaliszewski, A. Kumjian, J. Quigg and A. Sims, \emph{Topological realizations and fundamental groups of higher-rank graphs}, Proc. Edinb. Math. Soc., \textbf{59} (2016), 143--168.

\bibitem
{Kak} E.T.A. Kakariadis, \emph{Applications of entropy of product systems: higher-rank graphs}, Linear Algebra Appl., \textbf{594} (2020), 124--157.

\bibitem
{KV} J. Konter and A. Vdovina, \emph{Classifying polygonal algebras by their $K_0$-group}, Proc. Endin. Matrh. Soc. \textbf{58} (2015), 485--497.

\bibitem
{KP1} A. Kumjian and D. Pask, \emph{$C^*$-algebras of directed graphs and group actions}, Ergod. Th.  and Dynamical Systems \textbf{19 } (1999), 1503--1519.

\bibitem
{KP2} A. Kumjian and D. Pask, \emph{Higher rank graph $C^*$-algebras}, New York J. Math. \textbf{ 6} (2000), 1--20.

\bibitem
{KPRR} A. Kumjian, D. Pask,  I. Raeburn, and J. Renault, {\em Graphs, groupoids and Cuntz--Krieger algebras}, J. Funct. Anal. \textbf{ 144} (1997), 505--541.

\bibitem
{KPSW} A. Kumjian, D. Pask, A. Sims and M.F. Whittaker, \emph{Topological spaces associated to higher-rank graphs}. J. Comb. Th. Ser. A \textbf{143} (2016), 19--41.

\bibitem
{LV} M. Lawson and A. Vdovina. \emph{Higher dimensional generalisations of the Thompson groups}, Advances in Math. \textbf{ 369} (2020) 107--191.

\bibitem
{Ma} A. Ma{l}'cev, \emph{On the immersion of an algebraic ring into a field}, \newblock Math. Ann. \textbf{113} (1937), 686–-691.

\bibitem
{MRV} S.A. Mutter, A.-C. Radu and A. Vdovina, \emph{$C^*$-algebras of higher-rank graphs from groups acting on buildings, and explicit computation of their $K$-theory}, Publ. Mat. \textbf{68} (2024), 187--217.

\bibitem
{PQR1} D. Pask, I. Raeburn and J. Quigg, \emph{Fundamental Groupoids of $k$-graphs}, New York J. Math. \textbf{ 10} (2004), 195--207.

\bibitem
{PQR2} D. Pask, I. Raeburn and J. Quigg, \emph{Coverings of $k$-graphs}, J. Algebra \textbf{ 289} (2005), 161--191.

\bibitem
{PRW} D. Pask, I. Raeburn and N.A. Weaver, \emph{A family of $2$-graphs arising from two-dimensional subshifts}, Ergod. Th. and Dynamical Systems \textbf{ 29} (2009), 1613--1639.

\bibitem
{PRenS} D. Pask, A. Rennie and A. Sims, \emph{Noncommutative manifolds from graph and $k$-graph C*-algebras}, Comm. Math. Phys. \textbf{ 292} (2009) 607--636.

\bibitem
{RSY1} I. Raeburn, A. Sims  and T. Yeend, \emph{Higher-rank graphs and their {$C^*$}-algebras}, Proc. Edinb. Math. Soc. \textbf{46} (2003), 99--115.

\bibitem
{tfb} I. Raeburn and D.P. Williams, \emph{Morita equivalence and continuous-trace $C^*$-algebras}, Mathematical Surveys and Monographs, 60. American Mathematical Society, Providence, RI, 1998. xiv+327 pp. ISBN: 0--8218--0860--5.

\bibitem
{RobSteg1} G. Robertson and T. Steger, \emph{Affine buildings, tiling systems and higher rank Cuntz--Krieger algebras}, J. reine angew. Math.  \textbf{513} (1999), 115--144.

\bibitem
{RobSteg2} G. Robertson and T. Steger, \emph{Asymptotic $K$-theory for groups acting on $\tilde A_2$ buildings}, Can. J. Math. \textbf{ 53} (2001), 809--833.


\bibitem
{Ros} R. Rosjanuardi, \emph{Complex Kumjian--Pask algebras}, Acta Math. Sin. (Engl. Ser.) \textbf{29} (2013), 2073--2078.

\bibitem
{RSS} E. Ruiz, A. Sims and A.P.W. S\o rensen, \emph {U{CT}-{K}irchberg algebras have nuclear dimension one}, Adv. Math. \textbf{279} (2015), 1--28.

\bibitem
{Schubert} H. Schubert, Categories, Springer-Verlag, 1972.

\bibitem
{SZ} A. Skalski and J. Zacharias, \emph{Entropy of shifts on higher-rank graph $C^*$ -algebras}, Houston J. Math. \textbf{34} (2008),  269--282.

\bibitem
{Sp} J. Spielberg, \emph{Graph-based models for Kirchberg algebras}, J. Operator Theory \textbf{57} (2007), 347--374.

\bibitem
{V} A. Vdovina, \emph{Combinatorial structure of some hyperbolic buildings}, Math. Zeit. \textbf{241} (2002), 471--478.

\bibitem
{Y1} D. Yang, \emph{Endomorphisms and modular theory of $2$-graph $C^*$-algebras}, Indiana Univ. Math. J. \textbf{59} (2010) 495--520.

\bibitem
{Y4} D. Yang, \emph{The interplay between $k$-graphs and the Yang--Baxter equation}, J. Algebra \textbf{451} (2016), 494--525.

\end{thebibliography}
\end{document}